\documentclass[12pt, reqno]{amsart}
\usepackage{mathrsfs}
\usepackage{amssymb,amsthm,amsmath}
\usepackage[numbers,sort&compress]{natbib}
\usepackage{amssymb,amsmath}
\usepackage{amsfonts}
\usepackage{mathrsfs}
\usepackage{latexsym}
\usepackage{amssymb}
\usepackage{amsthm}
\usepackage{xcolor,booktabs,tabularx,array}
\usepackage{indentfirst}
\allowdisplaybreaks[4]
\renewcommand{\textwidth}{36pc}
\date{today}

\usepackage{hyperref}
\hypersetup{	colorlinks=true,
	linkcolor=blue,
	citecolor=red}

\usepackage{amsmath}
\usepackage{amsthm}

\theoremstyle{remark}
\newtheorem{remark}{Remark}[section]
\theoremstyle{plain}

\newtheorem{theorem}{Theorem}[section]
\newtheorem{proposition}{Proposition}[section]
\newtheorem{lemma}{Lemma}[section]
\newtheorem{corollary}{Corollary}[section]

\newcommand{\beq}{\begin{equation}}
	\newcommand{\eeq}{\end{equation}}
\newcommand{\ben}{\begin{eqnarray}}
	\newcommand{\een}{\end{eqnarray}}
\newcommand{\beno}{\begin{eqnarray*}}
	\newcommand{\eeno}{\end{eqnarray*}}

\numberwithin{equation}{section}

\begin{document}
	\title[Global solvability with arbitrary mass for a chemotaxis--fluid model]{Global solvability with arbitrary mass for a chemotaxis--fluid model with sensitivity growth between linear and quadratic}
	\author{Lili~Wang}
	\address[Lili~Wang]{State Key Laboratory of Mathematical Sciences, Academy of Mathematics $\&$ Systems Science, The Chinese Academy of Sciences, Beijing, 100190, China}
	\email{wllmath@163.com}
	\date{\today}
				\begin{abstract}
					We study a Patlak--Keller--Segel--Navier--Stokes system
					on $\mathbb T\times\mathbb R$ with density-dependent
					diffusion and chemotactic sensitivity growing like
					$n^{1+\alpha}$, where $0<\alpha<1$.
					The coefficients satisfy the constitutive relation
					of Painter and Hillen
					(Can. Appl. Math. Q. 2002; 10(4): 501--543).
					Near a sufficiently strong Couette flow, we prove
					global classical solvability and uniform density bounds
					for smooth initial data with arbitrary finite cell mass.
					A key ingredient is a joint energy estimate for the
					full density and the non-zero-mode vorticity, which
					preserves the positive nonlinear diffusion and avoids
					the mixed diffusion terms arising from separate
					estimates of the density modes.
					Bounded Fourier corrections yield enhanced dissipation
					within this energy, allowing the non-zero-mode
					density--fluid interactions to be controlled.
					For the density zero mode, which is not mixed by
					Couette transport, we observe that $\alpha<1$ allows
					one-dimensional interpolation to absorb its chemotactic
					production into diffusion for every finite mass.
					Together, these estimates control the nonlinear
					diffusion and fluid coupling without imposing
					a smallness condition on the initial cell mass.
				\end{abstract}
		\maketitle
		\color{black}
	\section{Introduction}
	Chemotaxis is the directed movement of cells in response to
	chemical signals and plays an important role in bacterial
	migration and aggregation. For example, Budrene and Berg \cite{BudreneBerg} observed that \emph{Escherichia coli} form spatial aggregates by moving toward a chemical attractant excreted by the cells themselves.
	In the Patlak--Keller--Segel model, diffusion competes with
	this self-attraction, which may lead to concentration and
	finite-time blow-up.
	The dependence of cell movement on local population density
	provides further motivation for nonlinear diffusion and
	chemotactic sensitivity.
	Painter and Hillen \cite{PainterHillen} developed such models
	from mechanisms involving finite cell size and cell density
	sensing.
	For cells suspended in a liquid, chemical cues also interact
	with fluid transport and density-induced buoyancy.
	Tuval et al. \cite{Tuval} studied suspensions of
	\emph{Bacillus subtilis} near the contact line of a sessile
	drop and showed that oxygen-directed swimming and gravity
	generate a persistent vortex that affects cell accumulation
	and enhances oxygen uptake.
	Their oxygen-consumption model illustrates the feedback
	between cell migration, chemical transport and fluid motion.
	This interaction motivates the mathematical study of how
	fluid transport modifies aggregation and whether sufficiently
	strong flows can prevent chemotactic blow-up.
	
For generalized chemotaxis models without fluid transport,
	the balance between nonlinear diffusion and chemotactic
	sensitivity plays an important role in aggregation.
	Ding and Winkler \cite{DW2024} studied a related
	parabolic--elliptic system in balls with Neumann boundary
	conditions.
	In two dimensions, when the ratio of sensitivity to diffusion
	is bounded below by a positive multiple of the density at
	high densities, sufficiently large masses admit radial
	initial data whose solutions become unbounded in finite
	or infinite time.
	Their chemical equation is a Poisson equation with the
	spatial average of the density subtracted from its source.
	Further finite-time blow-up results for chemotaxis systems
	without fluid transport, including parabolic--elliptic and
	quasilinear parabolic--parabolic models, can be found in
	\cite{Nagai2001,CS2014}.
	These results show that diffusion need not prevent
	chemotactic aggregation and motivate the study of the
	stabilizing effect of fluid transport.
	
	In two dimensions, blow-up suppression has been studied
both for prescribed velocity fields and for coupled
chemotaxis--fluid systems.
For prescribed velocities, Kiselev and Xu \cite{KX}
proved that sufficiently strong mixing flows prevent
chemotactic blow-up on $\mathbb T^2$ for arbitrary
nonnegative smooth initial densities.
Bedrossian and He \cite{BH} obtained global existence
without a mass restriction for a class of sufficiently
strong shear flows on $\mathbb T^2$.
Their argument combines enhanced dissipation of the
non-zero modes with estimates for the zero mode,
which is not mixed by the shear.
These works show that transport can accelerate
dissipation and suppress chemotactic aggregation.
For the two-dimensional
Patlak--Keller--Segel--Navier--Stokes system,
the velocity evolves with the cell density, so the
fluid perturbation and its feedback on chemotaxis
must also be controlled.
Zeng, Zhang and Zi \cite{ZZZ}
proved global regularity near sufficiently strong
Couette flow on $\mathbb T\times\mathbb R$
for the parabolic--elliptic system with constant
diffusion and linear sensitivity, without a smallness
condition on the initial cell mass.
Wang, Wang and Zhang \cite{WWZ} considered constant
diffusion and sensitivity $n^2+n$ on the same domain
and obtained global solutions under the sufficient
condition that the initial cell mass is less than
$\frac{2\pi}{\sqrt{3}}$.
In the whole plane, Chen, Wang, Yang and Zhang
\cite{CWYZ} obtained global bounded solutions with
arbitrary mass for a class of generalized
Hillen--Painter systems near Couette flow,
under anisotropic weighted assumptions on the initial
data. The domain and coefficient structure in that
result differ from those considered here. For further results on blow-up suppression in two dimensions,
we refer to \cite{He2018,HT2019,IXZ2021} for prescribed flows
and to \cite{He2023,CW2024} for coupled chemotaxis--fluid
systems, and the references therein.

	Related blow-up suppression results have also been
	established in three dimensions, both for prescribed
	shear flows \cite{BH,BHerratum} and for the fully
	coupled Patlak--Keller--Segel--Navier--Stokes system
	near Couette flow \cite{CWW3D}.
	These results impose mass conditions to control
	two-dimensional aggregation in the density zero mode.
	The coupled result also imposes smallness conditions
	on the initial velocity perturbation and certain
	density modes.
	Further results on blow-up suppression in three
	dimensions can be found in
	\cite{FSW2022,HTZ2022,He2025}
	and the references therein.
	
	These developments motivate the question addressed in this
	paper: {\bf can sufficiently strong Couette flow prevent blow-up
	on $\mathbb T\times\mathbb R$ for sensitivity growth strictly
	between the linear and quadratic cases, without restricting
	the initial cell mass?}
	We study this question with diffusion and sensitivity linked
	by the constitutive relation derived by Painter and Hillen
	\cite{PainterHillen}.
	To specify the model, we begin with the parabolic--elliptic
	Patlak--Keller--Segel--Navier--Stokes system
	\begin{equation}\label{eq:general-model}
		\left\{
		\begin{aligned}
			&\partial_t n+u\cdot\nabla n
			=\nabla\cdot\bigl(\phi(n)\nabla n\bigr)
			-\nabla\cdot\bigl(\psi(n)\nabla c\bigr),\\
			&\Delta c+n-c=0,\\
			&\partial_tu+u\cdot\nabla u+\nabla p
			=\Delta u+ne_2,\quad \nabla\cdot u=0,\\
			&(n,u)|_{t=0}=(n_{\mathrm{in}},u_{\mathrm{in}}),
			\quad n_{\mathrm{in}}\ge0,
		\end{aligned}
		\right.
	\end{equation}
	posed on $\mathbb T\times\mathbb R$, where
$\mathbb T=\mathbb R/(2\pi\mathbb Z)$ and $e_2=(0,1)$.
Here $n$, $c$, $u$ and $p$ denote the cell density,
the chemoattractant concentration, the fluid velocity and
the pressure, respectively.
The functions $\phi(n)$ and $\psi(n)$ describe
density-dependent diffusion and chemotactic sensitivity.
They belong to $C^2([0,\infty))$ and satisfy
$\phi>0$ on $[0,\infty)$, $\psi(0)=0$ and
$\psi>0$ on $(0,\infty)$.
The chemical concentration is determined by the elliptic
equation and has no independently prescribed initial datum.

Painter and Hillen \cite{PainterHillen} derived diffusion
and sensitivity coefficients from density-dependent cell
movement. We retain their constitutive relation
	\begin{equation}\label{eq:constitutive-revised}
		\phi(n)=Q(n)-nQ'(n),
		\qquad
		\psi(n)=nQ(n),
	\end{equation}
and consider the generalized choice
	\begin{equation*}
		Q(n)=(1+n)^\alpha,
		\qquad 0<\alpha<1.
	\end{equation*}
The sensitivity then grows like $n^{1+\alpha}$,
with total growth exponent $1+\alpha\in(1,2)$.
Substitution into \eqref{eq:constitutive-revised} gives
	\begin{equation}\label{eq:coefficients-revised}
		\phi(n)
		=(1+n)^{\alpha-1}\bigl[1+(1-\alpha)n\bigr],
		\qquad
		\psi(n)=n(1+n)^\alpha.
	\end{equation}
	Thus the constitutive relation determines both coefficients.
	For $0<\alpha<1$, the diffusion is positive and grows like
	$(1-\alpha)n^\alpha$ at high density.
	The change in sensitivity is therefore accompanied by
	a change in diffusion.
	With these coefficients, \eqref{eq:general-model} becomes
	\begin{equation}\label{eq:full-model-revised}
		\left\{
		\begin{aligned}
			&\partial_t n+u\cdot\nabla n=\nabla\cdot\!\left(
			(1+n)^{\alpha-1}
			\bigl[1+(1-\alpha)n\bigr]\nabla n
			\right)-\nabla\cdot\!\left(
			n(1+n)^\alpha\nabla c
			\right),\\
			&\Delta c+n-c=0,\\
			&\partial_tu+u\cdot\nabla u+\nabla p
			=\Delta u+ne_2,\quad \nabla\cdot u=0,\\
			&(n,u)|_{t=0}=(n_{\mathrm{in}},u_{\mathrm{in}}),
			\quad n_{\mathrm{in}}\ge0.
		\end{aligned}
		\right.
	\end{equation}
We study \eqref{eq:full-model-revised} near the Couette flow
$(Ay,0)$, where $A\ge1$, and write $u=(Ay,0)+v$.
To eliminate the pressure, we introduce the perturbation
vorticity $\omega=\partial_yv_1-\partial_xv_2$.
Let $\Phi$ be the stream function of $v$, so that
\[
v=\nabla^\perp\Phi=(\partial_{y}\Phi, -\partial_{x}\Phi),\qquad
\Delta\Phi=\omega.
\]
After the time rescaling $t\mapsto \frac{t}{A}$, with the rescaled
fields denoted by the same symbols, we obtain
\begin{equation}\label{eq:rescaled-system}
	\left\{
	\begin{aligned}
		&\partial_tn+y\partial_xn
		=\frac1A\nabla\cdot\!\left(
		(1+n)^{\alpha-1}
		\bigl[1+(1-\alpha)n\bigr]\nabla n
		\right)\\&\qquad\qquad\qquad-\frac1A\nabla\cdot
		\bigl(n(1+n)^\alpha\nabla c+vn\bigr),\\
		&\Delta c+n-c=0,\\
		&\partial_t\omega+y\partial_x\omega
		-\frac1A\Delta\omega=-\frac1A\partial_xn
		-\frac1A\nabla\cdot(v\omega),\\
		&(n,\omega)|_{t=0}=(n_{\mathrm{in}},\omega_{\mathrm{in}}).
	\end{aligned}
	\right.
\end{equation}
Here $0<\alpha<1$,
$v_{\mathrm{in}}=u_{\mathrm{in}}-(Ay,0)$ and
$\omega_{\mathrm{in}}
=\partial_yv_{1,\mathrm{in}}-\partial_xv_{2,\mathrm{in}}$.
	
	Before stating the result, we write
	$M=\|n_{\mathrm{in}}\|_{L^1(\mathbb T\times\mathbb R)}$. Moreover, for a given function $f=f(t,x,y)$, define its zero and non-zero modes by
	\begin{equation}\label{define f0}
			P_0f=f_0=\frac1{2\pi}\int_{\mathbb T}f(x,y)dx,
		\qquad P_{\neq}f=f_{\neq}=f-f_0.
	\end{equation}
	 Our objective is to prove global classical solvability and
	uniform density bounds for this coefficient pair when $A$
	is sufficiently large, with no smallness requirement on $M$
	or on the fixed initial perturbation.
	Before stating the global result, we record the local
	existence and continuation property used in its proof.
	This result holds for every fixed $A\ge1$; its proof is
	given in Appendix \ref{sec:local-continuation}.
	\begin{proposition}[Local existence and continuation]
		\label{prop:H3-continuation}
		Let $0<\alpha<1$ and $A\ge1$.
		Suppose that $0\le n_{\mathrm{in}}\in
		L^1\cap H^3(\mathbb T\times\mathbb R)$,
		$v_{\mathrm{in}}\in H^3(\mathbb T\times\mathbb R)$,
		and $\nabla\cdot v_{\mathrm{in}}=0$.
		Then system \eqref{eq:rescaled-system} has a unique
		local solution satisfying
		\[
		n,v\in C([0,T_*];H^3)\cap L^2(0,T_*;H^4)
		\]
		for some $T_*>0$. The solution is classical for positive times.
		If its maximal existence time $T_{\max}$ is finite and
		\begin{equation}\label{eq:continuation-density-assumption}
			\sup_{0\le t<T_{\max}}
			\bigl(\|n(t)\|_{L^\infty}+\|n(t)\|_{L^2}\bigr)<\infty,
		\end{equation}
		then it extends beyond $T_{\max}$.
	\end{proposition}
	The following theorem establishes this arbitrary-mass
	global result.
	\begin{theorem}\label{thm:full-global}
		Let $0<\alpha<1$. Suppose $0\leq n_{\mathrm{in}}\in L^1\cap H^3(\mathbb T\times\mathbb R), v_{\mathrm{in}}\in H^3(\mathbb T\times\mathbb R)$ and $\nabla\cdot v_{\mathrm{in}}=0$.
		There exists $A_0\ge1$, depending only on $\alpha$, $M$ and
		these initial norms, such that, for $A\ge A_0$,
		\eqref{eq:full-model-revised} with
		$u_{\mathrm{in}}=(Ay,0)+v_{\mathrm{in}}$ has a unique global
		classical solution. For every finite $T$, its perturbation
		satisfies $n,v\in C([0,T];H^3)\cap L^2(0,T;H^4)$.
		Moreover,
		\begin{align*}
			\sup_{t\ge0}\bigl(\|n(t)\|_{L^2}+\|n(t)\|_{L^\infty}
			+\|\omega_{\neq}(t)\|_{L^2}
			+\|\omega_0(t)\|_{L^\infty(\mathbb R)}\bigr)\le C\bigl(\alpha,M,\|n_{\mathrm{in}}\|_{H^3},
			\|v_{\mathrm{in}}\|_{H^3}\bigr).
		\end{align*}
		The density is nonnegative and preserves its initial mass.
	\end{theorem}
		\begin{remark}
		The nonlinear diffusion introduces a difficulty
		in estimating $n_{\neq}$ separately.
		Testing its equation produces a mixed term involving
		$\partial_yn_0\,\partial_yn_{\neq}$.
		Since its coefficient depends on the full density,
		this term does not in general vanish after integration
		in $x$ and has no definite sign.
		We therefore estimate the full density together with
		$\omega_{\neq}$ in the joint energy
		\eqref{eq:full-energy}.
		This retains the positive nonlinear diffusion integral
		and controls the density--velocity coupling.
		The density Fourier correction is divided by an
		upper bound for the diffusion coefficient.
		Its diffusion contribution can then be absorbed
		into the same positive integral without
		differentiating the coefficient.
	\end{remark}
	\begin{remark}
		The method also applies at $\alpha=0$.
		In this case the diffusion coefficient equals one
		and the chemotactic term is
		$-\nabla\cdot(n\nabla c)$.
		The zero-mode production satisfies
		\[
		\int_{\mathbb R}n_0^3dy
		\le C(M)
		\bigl(\|\partial_yn_0\|_{L^2(\mathbb R)}^2\bigr)^{\frac23}
		\le\varepsilon\|\partial_yn_0\|_{L^2(\mathbb R)}^2
		+C(M,\varepsilon).
		\]
		The joint energy and the zero-mode vorticity estimate
		remain valid.
		The gradient term used in the density iteration
		retains a positive coefficient, and the iteration
		exponents become $r_j=2^{j+1}$.
		Thus the argument also covers the classical linear
		sensitivity studied in \cite{ZZZ}.
		Theorem \ref{thm:full-global} is stated for
		$0<\alpha<1$ to focus on the intermediate
		superlinear sensitivities.
	\end{remark}
	\begin{remark}
		At $\alpha=1$, the diffusion coefficient equals one
		and the sensitivity is $n^2+n$, as in \cite{WWZ}.
		The distinction from $0<\alpha<1$ appears in
		two estimates.
		First, the zero-mode production contains a quartic
		term satisfying
		\[
		\int_{\mathbb R}n_0^4 dy
		\le C\left(\frac M{2\pi}\right)^2
		\|\partial_yn_0\|_{L^2(\mathbb R)}^2.
		\]
		The gradient exponent is now one, with a
		mass-dependent coefficient.
		Consequently, the arbitrary-mass absorption estimate
		\eqref{eq:all-mass-absorption} in
		Lemma \ref{lem:mean-production}
		does not extend to this endpoint.
		Increasing $A$ does not improve this absorption,
		since both production and diffusion carry
		the factor $A^{-1}$.
		Second, the higher-power gradient term in
		\eqref{eq:full-density-Moser-energy},
		used in Lemma \ref{lem:full-density-Moser},
		has a coefficient proportional to $1-\alpha$.
		At $\alpha=1$, this term disappears, while
		the same Young estimate produces
		$n^{r+2}|\nabla c|^2$.
		For $0<\alpha<1$, arbitrary-mass absorption
		and the positive higher-power gradient term
		are both available.
		These are the two estimates that allow the
		global argument to close without restricting
		the initial cell mass.
	\end{remark}
The following table compares the diffusion coefficients,
chemotactic sensitivities and mass conditions for
$\alpha=0$, $0<\alpha<1$ and $\alpha=1$.
The diffusion column gives the leading-order behavior
as $n\to\infty$.
The gradient exponent denotes the power of
$\|\partial_yn_0\|_{L^2(\mathbb R)}^2$
in the one-dimensional interpolation estimate for the
highest zero-mode chemotactic production term.
\begin{center}
	\small
	\renewcommand{\arraystretch}{1.5}
	\setlength{\tabcolsep}{4pt}
	\begin{tabularx}{\textwidth}{
			@{}cccc>{\centering\arraybackslash}X@{}
		}
		\toprule
		$\alpha$
		& \shortstack{Diffusion coefficient\\as $n\to\infty$}
		& \shortstack{Chemotactic\\sensitivity}
		& \shortstack{Zero-mode\\gradient exponent}
		& \shortstack{Mass condition\\in the stated results}\\
		\midrule
		
		$0$
		& $1$
		& $n$
		& $\frac23<1$
		& \mbox{No restriction \cite{ZZZ}}\\
		
		$(0,1)$
		& $(1-\alpha)n^\alpha$
		& $n(1+n)^\alpha$
		& $\frac{2+\alpha}{3}<1$
		& \mbox{No restriction, Theorem \ref{thm:full-global}}\\
		
		$1$
		& $1$
		& $n^2+n$
		& $1$
		& $M<\frac{2\pi}{\sqrt3}$ \cite{WWZ}\\
		
		\bottomrule
	\end{tabularx}
\end{center}

	\noindent\textbf{Notation.}
	\begin{itemize}
		\item For a given function $f=f(t,x,y)$, its zero mode $f_{0}$ and non-zero mode $f_{\neq}$ are defined by \eqref{define f0}.
		For the velocity components $v_1$ and $v_2$, we use
		$v_{1,0},v_{2,0}$ and $v_{1,\neq},v_{2,\neq}$ for their
		zero and non-zero modes. Since the vorticity is scalar,
		its modes are denoted by $\omega_0$ and $\omega_{\neq}$.
		\item At each fixed time, we use the Fourier transform
		\[
		\widehat f(t,k,\eta)
		=\frac1{2\pi}
		\int_{\mathbb R}\int_{\mathbb T}
		f(t,x,y){\rm e}^{-\mathrm i(kx+\eta y)}dxdy,
		\qquad k\in\mathbb Z,\quad \eta\in\mathbb R.
		\]
    	\item Zero-mode norms are taken on $\mathbb R$, and non-zero-mode norms on $\mathbb T\times\mathbb R$. 
		The interval in a space--time norm is the interval under consideration.
		\item The letter $C$ denotes a positive constant that may change from one line to the next. Unless additional arguments are displayed, it may
		depend on the fixed parameter $\alpha$ and, in $L^r$ estimates, on the fixed exponent $r$, but is independent of $A$ and the time interval. Dependence on $M$ or initial norms is displayed explicitly.
	\end{itemize}
	
	\noindent\textbf{Organization of the paper.}
	Section \ref{sec:main-ideas} introduces the energy and
	bootstrap bounds and explains the main ideas.
	Section \ref{sec:preliminaries} collects the coefficient,
	elliptic and zero-mode estimates.
	Section \ref{sec:full-energy} proves the joint
	density--vorticity energy estimate.
	Section \ref{sec:full-global-proof} establishes the
	uniform density bound and completes the proof of
	global existence. The local existence and continuation result is proved in Appendix \ref{sec:local-continuation}.
	\section{Main ideas and bootstrap bounds}
	\label{sec:main-ideas}
	
	The proof combines enhanced dissipation for the non-zero
	modes with one-dimensional estimates for the zero modes.
	The main difficulties are the coupling of density modes
	through nonlinear diffusion and the control of zero-mode
	chemotactic production for arbitrary mass.
	We first introduce the energy and bootstrap bounds,
	and then explain how these difficulties are addressed.
	\subsection{Energy and bootstrap bounds}
	Let $K_{\neq},K_\infty\ge1$ be fixed constants.
	For the maximal local solution supplied by
	Proposition\ref{prop:H3-continuation}, define
	\begin{equation}\label{eq:full-energy}
		\begin{aligned}
			E(t)={}&
			\|n(t)\|_{L^2}^2
			+\|\omega_{\neq}(t)\|_{L^2}^2+\sum_{k\ne0}\int_{\mathbb R}
			\operatorname{sgn}(k)
			\arctan\bigl(A^{-\frac13}|k|^{-\frac13}\eta\bigr)\\
			&\times
			\left(
			\frac{|\widehat n(t,k,\eta)|^2}{
				8\pi(1+2K_\infty)^{\alpha-1}
				[1+2(1-\alpha)K_\infty]}
			+\frac{|\widehat\omega(t,k,\eta)|^2}{8\pi}
			\right)d\eta.
		\end{aligned}
	\end{equation}
	The correction coefficients are constant in time.
	Since the density denominator is at least $8\pi$,
	Parseval's identity and
	$|\arctan z|\le\frac{\pi}{2}$ give
	\[
	\left|
	E(t)-\|n(t)\|_{L^2}^2-\|\omega_{\neq}(t)\|_{L^2}^2
	\right|
	\le
	\frac1{16}
	\left(
	\|n(t)\|_{L^2}^2+\|\omega_{\neq}(t)\|_{L^2}^2
	\right).
	\]
	Consequently,
	\begin{equation}\label{eq:full-energy-equivalence}
		\begin{aligned}
			\frac34
			\left(
			\|n(t)\|_{L^2}^2+\|\omega_{\neq}(t)\|_{L^2}^2
			\right)\le E(t)\le
			\frac54
			\left(
			\|n(t)\|_{L^2}^2+\|\omega_{\neq}(t)\|_{L^2}^2
			\right).
		\end{aligned}
	\end{equation}
	On a compact interval $[0,T]$, with $T<T_{\max}$,
	we assume the bootstrap bounds
	\begin{equation}\label{eq:full-bootstrap}
		E(t)\le2K_{\neq},
		\qquad
		\|n(t)\|_{L^\infty}\le2K_\infty,
		\qquad 0\le t\le T.
	\end{equation}
	All subsequent bootstrap estimates refer to
	\eqref{eq:full-bootstrap}.
	The constants $K_{\neq}$ and $K_\infty$ will be
	chosen from the initial data before fixing the
	lower bound for $A$.
	\subsection{A joint estimate for density and vorticity}
	The density-dependent diffusion contains both
	$n_0$ and $n_{\neq}$.
	Estimating the density modes separately therefore
	introduces mixed diffusion terms involving
	$\partial_y n_0$ and $\nabla n_{\neq}$,
	whose contributions have no definite sign.
	
	\noindent\textbf{Key observation 1: retaining the full
		density preserves the nonlinear diffusion.}
	We estimate the full density together with
	$\omega_{\neq}$ in the energy \eqref{eq:full-energy}.
	Testing the density equation with $2n$ retains
	the positive term
	\[
	\frac2A\int_{\mathbb T\times\mathbb R}
	(1+n)^{\alpha-1}[1+(1-\alpha)n]
	|\nabla n|^2 dxdy.
	\]
	Thus the mixed diffusion terms do not need to be
	estimated separately.
	The denominator in the density Fourier correction
	contains the upper bound of the diffusion coefficient
	on $0\le n\le2K_\infty$.
	Using the boundedness of the multiplier and integration
	by parts, its diffusion contribution is controlled
	by the displayed dissipation without differentiating
	the coefficient.
	The arctangent corrections also yield dissipation
	from Couette transport.
	Together with ordinary diffusion, they control
	the non-zero modes at the scale $A^{-\frac13}$.
	
	Including $\omega_{\neq}$ in the same energy
	controls the coupling to the fluid.
	For the velocity contribution in the density correction,
	the Biot--Savart estimate and
	\eqref{eq:full-bootstrap} give
	\[
	\|v_{\neq}n\|_{L^2}
	\le CK_\infty\|\omega_{\neq}\|_{L^2}.
	\]
	The density forcing in the vorticity equation is
	estimated by
	\[
	\frac CA
	\|n_{\neq}\|_{L^2}
	\|\nabla\omega_{\neq}\|_{L^2}
	\le
	\frac{\varepsilon}{A}
	\|\nabla\omega_{\neq}\|_{L^2}^2
	+\frac{C_\varepsilon}{A}
	\|n_{\neq}\|_{L^2}^2.
	\]
	The resulting gradient terms are absorbed by ordinary
	diffusion, while the remaining terms with coefficient
	$A^{-1}$ are absorbed by enhanced dissipation for
	sufficiently large $A$.
	The complete joint estimate is proved in
	Section \ref{sec:full-energy}.
	
	\subsection{Estimates for the zero modes}
	
	Couette transport does not act on the zero modes.
	The density zero mode must therefore be controlled
	by diffusion, while the vorticity zero mode supplies
	the coefficient bound needed to estimate the
	zero-mode velocity transport.
	
	\noindent\textbf{Key observation 2: $\alpha<1$ allows
		zero-mode absorption for arbitrary mass.}
	Using $\Delta c=c-n$ and $c\ge0$, the chemotactic
	contribution to the density energy is bounded by
	\[
	\frac{C_\alpha}{A}
	\int_{\mathbb T\times\mathbb R}
	\bigl(n^3+n^{3+\alpha}\bigr)dxdy.
	\]
The density bootstrap bound and the identity
$\int_{\mathbb T}n_{\neq}dx=0$ imply
	\[
	\begin{aligned}
		\int_{\mathbb T\times\mathbb R}
		\bigl(n^3+n^{3+\alpha}\bigr)dxdy\le
		2\pi\int_{\mathbb R}
		\bigl(n_0^3+n_0^{3+\alpha}\bigr)dy
		+C(\alpha,K_\infty)\|n_{\neq}\|_{L^2}^2.
	\end{aligned}
	\]
	The non-zero-mode remainder is controlled by
	enhanced dissipation.
	For the zero-mode contribution, mass conservation
	and one-dimensional interpolation yield
	\[
	\int_{\mathbb R}n_0^{3+\alpha}dy
	\le
	C(\alpha,M)
	\left(
	\|\partial_y n_0\|_{L^2(\mathbb R)}^2
	\right)^{\frac{2+\alpha}{3}}.
	\]
	Since $\frac{2+\alpha}{3}<1$, for every $\varepsilon>0$, Young's inequality gives
	\[
	C_\alpha\int_{\mathbb R}
	\bigl(n_0^3+n_0^{3+\alpha}\bigr)dy
	\le
	\varepsilon\|\partial_y n_0\|_{L^2(\mathbb R)}^2
	+C(\alpha,M,\varepsilon).
	\]
	This estimate absorbs the zero-mode production
	into density diffusion for every finite mass.
	It is the step that removes the mass restriction.
	The detailed interpolation estimates are given in
	Lemma \ref{lem:mean-production}.
	
	We also need a uniform bound for $\omega_0$.
	The contribution of
	$A^{-1}v_{1,0}\partial_x$ to the Fourier corrections
	is estimated by a commutator bound involving only
	$\|\partial_y v_{1,0}\|_{L^\infty(\mathbb R)}=
	\|\omega_0\|_{L^\infty(\mathbb R)}$.
	Thus a uniform bound for the amplitude of
	$v_{1,0}$ is unnecessary.
	Incompressibility gives
	\[
	\partial_t\omega_0-\frac1A\partial_y^2\omega_0
	=
	-\frac1A\partial_y^2(v_{1,\neq}v_{2,\neq})_0.
	\]
	The Biot--Savart estimates imply
	\[
	\begin{aligned}
		\|(v_{1,\neq}v_{2,\neq})_0\|_{L^1(\mathbb R)}
		\le C\|\omega_{\neq}\|_{L^2}^2,\quad
		\|\partial_y(v_{1,\neq}v_{2,\neq})_0
		\|_{L^{\frac43}(\mathbb R)}
		\le C\|\omega_{\neq}\|_{L^2}^2.
	\end{aligned}
	\]
	In Duhamel's formula, write $\tau=\frac{t-s}{A}$.
	For $0<\tau\le1$, placing one derivative on the
	averaged product gives the integrable heat-kernel
	factor $\tau^{-\frac78}$.
	For $\tau>1$, placing both derivatives on the
	heat kernel gives $\tau^{-\frac32}$, which is integrable
	at infinity.
	The change of variables cancels the factor $A^{-1}$
	in the forcing.
	Lemma \ref{lem:mean-vorticity-uniform} therefore gives
	\[
	\|\omega_0(t)\|_{L^\infty(\mathbb R)}
	\le
	\|(\omega_{\mathrm{in}})_0\|_{L^\infty(\mathbb R)}
	+C\sup_{0\le s\le t}
	\|\omega_{\neq}(s)\|_{L^2}^2,
	\]
	with $C$ independent of $A$ and time.
	This is the zero-mode bound required in the
	joint energy estimate.
	
	Combining the joint estimate with the zero-mode
	bounds yields
	\[
	E'(t)+\frac1A E(t)
	\le\frac{C(\alpha,M)}A
	\]
	on the bootstrap interval, provided that $A$ is
	sufficiently large in terms of the fixed bootstrap
	constants.
	The additive constant is independent of
	$K_{\neq}$ and $K_\infty$.
	The energy equivalence
	\eqref{eq:full-energy-equivalence} therefore allows
	$K_{\neq}$ to be chosen from the initial data
	independently of $K_\infty$ and $A$.
	
	The full-density Moser estimate in
	Lemma \ref{lem:full-density-Moser}
	then determines $K_\infty$ from this energy bound
	and the initial density.
	Finally, $A$ is chosen to satisfy the absorption
	conditions for these fixed constants.
	Section \ref{sec:full-global-proof} proves the
	strict bootstrap improvements and uses
	Proposition \ref{prop:H3-continuation}
	to complete the global existence proof.

\section{Preliminary estimates}
\label{sec:preliminaries}
We first establish the coefficient bounds and
conservation laws. We then collect the elliptic and
velocity estimates and derive the zero-mode bounds
needed in the joint energy argument.
\subsection{Coefficient bounds and conservation laws}
\begin{lemma}\label{lem:revised-coefficients}
	The functions
	\[
	s\longmapsto
	(1+s)^{\alpha-1}\bigl[1+(1-\alpha)s\bigr],
	\qquad
	s\longmapsto s(1+s)^\alpha
	\]
	are smooth on $[0,\infty)$.
	The diffusion coefficient is nondecreasing, and, for every
	$s\ge0$, there holds
		\begin{align*}
			&(1+s)^{\alpha-1}\bigl[1+(1-\alpha)s\bigr]\ge1,\\
			&(1+s)^{\alpha-1}\bigl[1+(1-\alpha)s\bigr]
			\ge\frac{1-\alpha}{2}(1+s^\alpha),\\
			&\frac{s^2(1+s)^{2\alpha}}
			{(1+s)^{\alpha-1}\bigl[1+(1-\alpha)s\bigr]}
			\le C_\alpha\bigl(s^2+s^{2+\alpha}\bigr).
		\end{align*}
\end{lemma}
\begin{proof}
	For $s\ge0$, direct calculations give
	\[
	\frac{d}{ds}
	\left\{
	(1+s)^{\alpha-1}\bigl[1+(1-\alpha)s\bigr]
	\right\}
	=\alpha(1-\alpha)s(1+s)^{\alpha-2}\ge0.
	\]
  Since the diffusion coefficient equals one at $s=0$,
	this proves the first bound.
	Note that
	\[
	1+(1-\alpha)s\ge(1-\alpha)(1+s),
	\qquad
	(1+s)^\alpha\ge\frac12(1+s^\alpha).
	\]
	These inequalities prove the second bound.
	Using the same lower bound for the diffusion coefficient
	and $(1+s)^\alpha\le1+s^\alpha$, we obtain
	\[
	\frac{s^2(1+s)^{2\alpha}}
	{(1+s)^{\alpha-1}\bigl[1+(1-\alpha)s\bigr]}
	\le
	\frac{s^2(1+s)^\alpha}{1-\alpha}
	\le
	\frac{s^2+s^{2+\alpha}}{1-\alpha}.
	\]
\end{proof}
The estimates for the density zero mode use nonnegativity
and mass conservation. Incompressibility also determines
the velocity zero modes.
\begin{lemma}\label{lem:revised-mass}
	The density is nonnegative and
	\begin{equation*}
		\int_{\mathbb T\times\mathbb R}
		n(t,x,y)dxdy=M,
		\qquad
		\int_{\mathbb R}n_0(t,y)dy=\frac{M}{2\pi}.
	\end{equation*}
	Moreover, $v_{2,0}=0$ and $\partial_yv_{1,0}=\omega_0$.
\end{lemma}
\begin{proof}
	Fix $0<T<T_{\max}$.
	By the local existence assertion of
	Proposition \ref{prop:H3-continuation},
	proved in Appendix \ref{sec:local-continuation},
	we have
	$n,v\in C([0,T];H^3(\mathbb T\times\mathbb R))$.
	The chemical equation gives
	\[
	\widehat c(t,k,\eta)
	=\frac{\widehat n(t,k,\eta)}{1+k^2+\eta^2},
	\qquad
	\|c(t)\|_{H^5}\le C\|n(t)\|_{H^3}.
	\]
	Consequently, Sobolev embedding yields
	\begin{equation}\label{eq:mass-regularity-bounds}
		\begin{aligned}
			\sup_{0\le t\le T}\bigl(
			&\|n(t)\|_{L^\infty}
			+\|v(t)\|_{L^\infty}
			+\|\nabla c(t)\|_{L^\infty}\\
			&+\|n(t)\|_{L^2}
			+\|\nabla n(t)\|_{L^2}
			+\|v(t)\|_{L^2}
			+\|\nabla c(t)\|_{L^2}
			\bigr)\le C_T.
		\end{aligned}
	\end{equation}
	The constants $C_T$ below may depend on the solution
	on $[0,T]$, but are independent of $R$.
	Choose an even function $\chi\in C_c^\infty(\mathbb R)$
	such that $0\le\chi\le1$, $\chi=1$ on $[-1,1]$,
	$\chi=0$ outside $[-2,2]$, and $\chi$ is nonincreasing
	on $[0,\infty)$.
	Then $\chi(\frac yR)$ increases pointwise to one as
	$R\to\infty$.
	Since
	$\partial_y\left[\chi\left(\frac yR\right)\right]
	=\frac1R\chi'\left(\frac yR\right)$,
	a change of variables gives
	\begin{equation}\label{eq:cutoff-derivative-bounds}
		\left\|
		\partial_y\left[\chi\left(\frac yR\right)\right]
		\right\|_{L^2(\mathbb R)}
		=R^{-\frac12}\|\chi'\|_{L^2(\mathbb R)},
		\quad
		\left\|
		\partial_y\left[\chi\left(\frac yR\right)\right]
		\right\|_{L^\infty(\mathbb R)}
		=R^{-1}\|\chi'\|_{L^\infty(\mathbb R)}.
	\end{equation}
	
	\noindent\textbf{Nonnegativity.}
	The continuity of $n$ in $H^3$, supplied by
	Proposition \ref{prop:H3-continuation},
	and Sobolev embedding give
	$$\|n(t)-n_{\mathrm{in}}\|_{L^\infty}
	\le C\|n(t)-n_{\mathrm{in}}\|_{H^3}
	\to0\quad\text{as}\quad t\to0^+.$$
	Due to $n_{\mathrm{in}}\ge0$, it follows that
	\[
	\inf_{(x,y)\in\mathbb T\times\mathbb R}n(t,x,y)
	\ge-\|n(t)-n_{\mathrm{in}}\|_{L^\infty}>-\frac12
	\]
	on a sufficiently short initial time interval.
	We first prove nonnegativity as long as this lower
	bound holds.
	For $-\frac12<n<0$, the explicit coefficients satisfy
	\begin{equation}\label{diffion coe}
		\frac12
		\le(1+n)^{\alpha-1}[1+(1-\alpha)n]
		\le2^{1-\alpha},
		\qquad
		|n(1+n)^\alpha|\le|\min\{n,0\}|.
	\end{equation}
	The bounds in \eqref{diffion coe} apply to all terms
	below containing $\min\{n,0\}$ or its gradient,
	since both vanish almost everywhere on $\{n\ge0\}$.
	
	Test $\eqref{eq:rescaled-system}_1$ with
	$\chi(\frac yR)\min\{n,0\}$.
	For almost every $t$, integration by parts,
	periodicity in $x$ and $\nabla\cdot v=0$ give
	\begin{equation}\label{energy ineq}
		\begin{aligned}
			&\frac12\frac d{dt}
			\int_{\mathbb T\times\mathbb R}
			\chi\left(\frac yR\right)
			|\min\{n,0\}|^2 dxdy\\
			&+
			\frac1A\int_{\mathbb T\times\mathbb R}
			\chi\left(\frac yR\right)
			(1+n)^{\alpha-1}[1+(1-\alpha)n]
			|\nabla\min\{n,0\}|^2 dxdy\\=~&
			\frac1A\int_{\mathbb T\times\mathbb R}
			\chi\left(\frac yR\right)n(1+n)^\alpha
			\nabla c\cdot\nabla\min\{n,0\}dxdy\\
			&-\frac1{A}\int_{\mathbb T\times\mathbb R}
			\partial_y\left[\chi\left(\frac yR\right)\right]
			(1+n)^{\alpha-1}[1+(1-\alpha)n]
			\partial_yn\min\{n,0\}dxdy\\
			&+\frac1{A}\int_{\mathbb T\times\mathbb R}
			\partial_y\left[\chi\left(\frac yR\right)\right]
			n(1+n)^\alpha
			\partial_yc\,\min\{n,0\}dxdy\\
			&+\frac1{2A}\int_{\mathbb T\times\mathbb R}
			\partial_y\left[\chi\left(\frac yR\right)\right]
			v_2|\min\{n,0\}|^2 dxdy
			=:I_{1}+I_{2}+I_{3}+I_{4}.
		\end{aligned}
	\end{equation}
	Using the bound for $n(1+n)^\alpha$ in
	\eqref{diffion coe}, the bound for $\nabla c$ in
	\eqref{eq:mass-regularity-bounds}
	and Young's inequality, we obtain
	\begin{align*}
		I_{1}\le~&
		\frac{\|\nabla c(t)\|_{L^\infty}}A
		\int_{\mathbb T\times\mathbb R}
		\chi\left(\frac yR\right)
		|\min\{n,0\}||\nabla\min\{n,0\}|dxdy\\\le~&
		\frac1{4A}\int_{\mathbb T\times\mathbb R}
		\chi\left(\frac yR\right)
		|\nabla\min\{n,0\}|^2 dxdy+
		\frac{C_T}{A}\int_{\mathbb T\times\mathbb R}
		\chi\left(\frac yR\right)
		|\min\{n,0\}|^2 dxdy.
	\end{align*}
	For the three terms containing the cutoff derivative,
	we first estimate each term in absolute value.
	The coefficient bounds in \eqref{diffion coe} give
	\begin{align*}
		&I_{2}+I_{3}+I_{4}
		\le |I_2|+|I_3|+|I_4|\\\le~&
		\frac CA
		\left\|
		\partial_y\left[\chi\left(\frac yR\right)\right]
		\right\|_{L^\infty(\mathbb R)}
		\int_{\mathbb T\times\mathbb R}
		\left[
		|\partial_yn|\,|\min\{n,0\}|
		+\bigl(|\partial_yc|+|v_2|\bigr)
		|\min\{n,0\}|^2
		\right]dxdy\\\le~&
		\frac{C_T}{AR}\int_{\mathbb T\times\mathbb R}
		\left(
		|\partial_yn|\,|\min\{n,0\}|
		+|\min\{n,0\}|^2
		\right)dxdy\\\le~&
		\frac{C_T}{AR}
		\left(
		\|\partial_yn(t)\|_{L^2}
		\|\min\{n(t),0\}\|_{L^2}
		+\|\min\{n(t),0\}\|_{L^2}^2
		\right)\\\le~&
		\frac{C_T}{AR}
		\left(
		\|\nabla n(t)\|_{L^2}\|n(t)\|_{L^2}
		+\|n(t)\|_{L^2}^2
		\right).
	\end{align*}
	Here the third line uses the $L^\infty$ derivative
	bound in \eqref{eq:cutoff-derivative-bounds}
	and the bounds for $\nabla c$ and $v$ in
	\eqref{eq:mass-regularity-bounds}.
	The fourth line follows from Cauchy--Schwarz,
	and the last uses $|\min\{n,0\}|\le|n|$.
	Substituting these estimates into \eqref{energy ineq}
	and using the diffusion lower bound in
	\eqref{diffion coe}, we retain the gradient term
	with coefficient $\frac1{4A}$.
	Integrating from $0$ to $t$ and using
	$n_{\mathrm{in}}\ge0$, we obtain
	\begin{equation}\label{eq:negative-integrated}
		\begin{aligned}
			&\frac12\int_{\mathbb T\times\mathbb R}
			\chi\left(\frac yR\right)
			|\min\{n(t,x,y),0\}|^2 dxdy\\&+
			\frac1{4A}\int_0^t\int_{\mathbb T\times\mathbb R}
			\chi\left(\frac yR\right)
			|\nabla\min\{n(s,x,y),0\}|^2 dxdyds\\\le~&
			\frac{C_T}{A}\int_0^t\int_{\mathbb T\times\mathbb R}
			\chi\left(\frac yR\right)
			|\min\{n(s,x,y),0\}|^2 dxdyds\\&+
			\frac{C_T}{AR}\int_0^t
			\left(
			\|n(s)\|_{L^2}\|\nabla n(s)\|_{L^2}
			+\|n(s)\|_{L^2}^2
			\right)ds.
		\end{aligned}
	\end{equation}
	By \eqref{eq:mass-regularity-bounds}, the final time
	integral is bounded uniformly for $0\le t\le T$.
	Hence its prefactor $R^{-1}$ makes the corresponding
	term tend to zero.
	Letting $R\to\infty$ in \eqref{eq:negative-integrated},
	using monotone convergence and discarding the
	nonnegative gradient term, we obtain
	\begin{equation*}
		\int_{\mathbb T\times\mathbb R}
		|\min\{n(t,x,y),0\}|^2dxdy
		\le
		\frac{C_T}{A}\int_0^t\int_{\mathbb T\times\mathbb R}
		|\min\{n(s,x,y),0\}|^2dxdyds.
	\end{equation*}
	Gronwall's inequality gives
	$\min\{n(t),0\}=0$ almost everywhere.
	Since $n$ is continuous in space, we conclude that
	$n\ge0$ on this interval.
	
	If the condition
	$\inf_{\mathbb T\times\mathbb R}n>-\frac12$
	had a first endpoint before $T$,
	continuity in $L^\infty$ would give $n\ge0$
	at that endpoint and preserve $n>-\frac12$
	on a further time interval.
	This contradicts the choice of the endpoint.
	Therefore $n\ge0$ throughout $[0,T]$.
	
	\noindent\textbf{Mass conservation.}
	Integrating the density equation
	$\eqref{eq:rescaled-system}_1$ against
	$\chi(\frac yR)$ and then integrating in time gives
	\begin{align*}
		&\int_{\mathbb T\times\mathbb R}
		n(t,x,y)\chi\left(\frac yR\right)dxdy
		-\int_{\mathbb T\times\mathbb R}
		n_{\mathrm{in}}(x,y)\chi\left(\frac yR\right)dxdy\\=&
		-\frac1A\int_0^t\int_{\mathbb T\times\mathbb R}
		\partial_y\left[\chi\left(\frac yR\right)\right]
		(1+n)^{\alpha-1}[1+(1-\alpha)n]
		\partial_yn dxdyds\\
		&+\frac1A\int_0^t\int_{\mathbb T\times\mathbb R}
		\partial_y\left[\chi\left(\frac yR\right)\right]
		\bigl[n(1+n)^\alpha\partial_yc+v_2n\bigr]dxdyds.
	\end{align*}
	Since $n$ is nonnegative and bounded by
	\eqref{eq:mass-regularity-bounds}, we have
	\[
	(1+n)^{\alpha-1}[1+(1-\alpha)n]\le C_T,
	\qquad
	n(1+n)^\alpha\le C_Tn.
	\]
	Noting that the cutoff is independent of $x$ and
	the length of $\mathbb T$ is $2\pi$, we have
	$\left\|
	\partial_y\left[\chi\left(\frac yR\right)\right]
	\right\|_{L^2(\mathbb T\times\mathbb R)}
	=(2\pi)^{\frac12}
	\left\|
	\partial_y\left[\chi\left(\frac yR\right)\right]
	\right\|_{L^2(\mathbb R)}$.
	Moreover, Cauchy--Schwarz gives
	\[
	\|n(1+n)^\alpha\partial_yc\|_{L^1}
	\le C_T\|n\|_{L^2}\|\nabla c\|_{L^2},
	\qquad
	\|v_2n\|_{L^1}
	\le\|v\|_{L^2}\|n\|_{L^2}.
	\]
	Combining these estimates and applying
	\eqref{eq:cutoff-derivative-bounds}, we obtain
	\begin{align*}
		&\left|
		\int_{\mathbb T\times\mathbb R}
		n(t,x,y)\chi\left(\frac yR\right)dxdy
		-
		\int_{\mathbb T\times\mathbb R}
		n_{\mathrm{in}}(x,y)\chi\left(\frac yR\right)dxdy
		\right|\\\le~&
		\frac{C_T}{A}
		\left\|
		\partial_y\left[\chi\left(\frac yR\right)\right]
		\right\|_{L^2(\mathbb R)}
		\int_0^t\|\nabla n(s)\|_{L^2}ds\\&+
		\frac{C_T}{A}
		\left\|
		\partial_y\left[\chi\left(\frac yR\right)\right]
		\right\|_{L^\infty(\mathbb R)}
		\int_0^t
		\left(
		\|n(s)\|_{L^2}\|\nabla c(s)\|_{L^2}
		+\|v(s)\|_{L^2}\|n(s)\|_{L^2}
		\right)ds\\\le~&
		\frac{C_T}{A}R^{-\frac12}
		\int_0^t\|\nabla n(s)\|_{L^2}ds+
		\frac{C_T}{A}R^{-1}
		\int_0^t
		\left(
		\|n(s)\|_{L^2}\|\nabla c(s)\|_{L^2}
		+\|v(s)\|_{L^2}\|n(s)\|_{L^2}
		\right)ds\\\le~&
		\frac{C_T}{A}(R^{-\frac12}+R^{-1}).
	\end{align*}
	The last inequality follows from the $L^2$ bounds in
	\eqref{eq:mass-regularity-bounds} and $t\le T$.
	Let $R\to\infty$.
	The initial cutoff integral tends to $M$,
	and the right-hand side tends to zero.
	Since $n\ge0$ and $\chi(\frac yR)\to1$,
	monotone convergence gives
	$\int_{\mathbb T\times\mathbb R}
	n(t,x,y)dxdy=M<\infty$.
	In particular, $n(t)\in L^1$. Therefore
	\[
	\int_{\mathbb R}n_0(t,y)dy
	=\frac1{2\pi}\int_{\mathbb R}\int_{\mathbb T}
	n(t,x,y)dxdy
	=\frac{M}{2\pi}.
	\]
	
	\noindent\textbf{Velocity zero modes.}
	Incompressibility and periodicity give
	\[
	\partial_yv_{2,0}
	=-\frac1{2\pi}\int_{\mathbb T}
	\partial_xv_1(t,x,y)dx=0.
	\]
	Moreover, Jensen's inequality yields
	$\|v_{2,0}(t)\|_{L^2(\mathbb R)}
	\le(2\pi)^{-\frac12}
	\|v_2(t)\|_{L^2(\mathbb T\times\mathbb R)}$.
	Thus $v_{2,0}$ is constant in $y$ and belongs to
	$L^2(\mathbb R)$, which implies $v_{2,0}=0$.
	Finally, averaging the definition of vorticity gives
	\[
	\omega_0
	=
	\partial_yv_{1,0}
	-\frac1{2\pi}\int_{\mathbb T}
	\partial_xv_2(t,x,y)dx
	=
	\partial_yv_{1,0}.
	\]
\end{proof}
\subsection{Elliptic and velocity estimates}
The chemical equation is unchanged by the nonlinear
coefficient choice. For the non-zero modes, its Fourier
symbol and Parseval's identity give
\begin{equation}\label{eq:nonzero-chemical-gradient}
	\begin{aligned}
		\|\nabla c_{\neq}(t)\|_{L^2}^2=\sum_{k\ne0}\int_{\mathbb R}
		\frac{k^2+\eta^2}{(1+k^2+\eta^2)^2}
		|\widehat n(t,k,\eta)|^2 d\eta\le\frac14\|n_{\neq}(t)\|_{L^2}^2.
	\end{aligned}
\end{equation}
Noting that
\[
\|\nabla c\|_{L^4}\le
(2\pi)^{\frac14}\|\partial_yc_0\|_{L^4(\mathbb R)}
+\|\nabla c_{\neq}\|_{L^4},
\]
the elliptic estimates in
\cite[Lemmas 3.2--3.3]{WWZ} imply
\begin{equation}\label{eq:full-elliptic-gradient}
	\|\nabla c\|_{L^4}\le C\|n\|_{L^2}.
\end{equation}
For the velocity, the Biot--Savart calculation in
\cite[proof of Lemma 3.5, equations (3.2)--(3.5)]{WWZ}
and the embedding $H^1\hookrightarrow L^4$ give
\begin{equation}\label{eq:Biot-bound}
	\|v_{\neq}\|_{L^2}+\|v_{\neq}\|_{L^4}
	\le C\|\omega_{\neq}\|_{L^2}.
\end{equation}
We also use the interpolation inequality on
$\mathbb T\times\mathbb R$:
\begin{equation}\label{eq:cylinder-interpolation}
	\|f\|_{L^4}^2
	\le C\|f\|_{L^2}
	\bigl(\|\nabla f\|_{L^2}+\|f\|_{L^2}\bigr),
	\qquad f\in H^1(\mathbb T\times\mathbb R).
\end{equation}
The additional $L^2$ term accounts for the periodic direction.
\begin{lemma}\label{lem:mass-elliptic}
	The chemical concentration is nonnegative and
	\begin{equation}\label{eq:mass-chemical}
		\|\partial_yc_0(t)\|_{L^\infty(\mathbb R)}
		\le\frac{M}{4\pi}.
	\end{equation}
	Moreover,
	\[
	\|c(t)\|_{H^{j+2}}\le C_j\|n(t)\|_{H^j},
	\qquad j\ge0,
	\]
	whenever $n(t)\in H^j$, and
	\[
	\|\nabla c(t)\|_{L^\infty}
	\le C\|n(t)\|_{L^\infty}.
	\]
\end{lemma}

\begin{proof}
	Fix a time $t$. The elliptic equation gives
	$\widehat c(t,k,\eta)=\frac{\widehat n(t,k,\eta)}{1+k^2+\eta^2}$.
	Parseval's identity proves the Sobolev estimate.
	The same formula gives the representation
	\begin{equation}\label{heat}
		c(t)=\int_0^\infty {\rm e}^{-s}{\rm e}^{s\Delta}n(t)ds
	\end{equation}
	in $L^2$. Since $n(t)\ge0$ and the heat kernel is nonnegative,
	we obtain $c(t)\ge0$.
	
	The one-dimensional Green function gives
	$c_0(t,y)=\frac12\int_{\mathbb R}
	{\rm e}^{-|y-z|}n_0(t,z)dz$.
	Its derivative is
	\[
	\partial_yc_0(t,y)
	=-\frac12\int_{\mathbb R}
	\operatorname{sgn}(y-z){\rm e}^{-|y-z|}
	n_0(t,z)dz.
	\]
	Therefore
	\[
	|\partial_yc_0(t,y)|
	\le\frac12\int_{\mathbb R}n_0(t,z)dz
	=\frac{M}{4\pi},
	\]
	which proves \eqref{eq:mass-chemical}.
	The heat kernel on $\mathbb T\times\mathbb R$ is
	\[
	\sum_{j\in\mathbb Z}
	\frac1{4\pi s}
	\exp\left(
	-\frac{(x+2\pi j)^2+y^2}{4s}
	\right),
	\qquad s>0.
	\]
	By the triangle inequality, the $L^1$ norm of its
	gradient on $\mathbb T\times\mathbb R$ is bounded by
	the $L^1$ norm of the gradient of the Gaussian on
	$\mathbb R^2$. Scaling gives
	\[
	\int_{\mathbb R^2}
	\left|
	\nabla\left[
	\frac1{4\pi s}
	\exp\left(-\frac{x^2+y^2}{4s}\right)
	\right]
	\right|dxdy
	\le Cs^{-\frac12}.
	\]
	Since the heat semigroup acts by convolution with
	this kernel, Young's convolution inequality yields
	\[
	\|\nabla e^{s\Delta}f\|_{L^\infty}
	\le Cs^{-\frac12}\|f\|_{L^\infty},
	\qquad s>0.
	\]
	Recall the heat representation in \eqref{heat}.
	The preceding estimate gives
	\[
	e^{-s}\|\nabla {\rm e}^{s\Delta}n(t)\|_{L^\infty}
	\le C{\rm e}^{-s}s^{-\frac12}\|n(t)\|_{L^\infty}.
	\]
	We therefore differentiate the heat representation 
	under the integral sign to obtain
	$\nabla c(t)
	=\int_0^\infty {\rm e}^{-s}\nabla {\rm e}^{s\Delta}n(t)ds$.
	Consequently,
	\[
	\|\nabla c(t)\|_{L^\infty}
	\le
	C\|n(t)\|_{L^\infty}
	\int_0^\infty {\rm e}^{-s}s^{-\frac12}ds
	\le C\|n(t)\|_{L^\infty}.
	\]
\end{proof}
\subsection{Estimates for the density zero mode}
\begin{lemma}\label{lem:mean-production}
There holds
	\begin{equation}\label{eq:ordinary-mean-GN}
		\|n_0\|_{L^2(\mathbb R)}^6
		\le
		C\left(\frac{M}{2\pi}\right)^4
		\|\partial_yn_0\|_{L^2(\mathbb R)}^2.
	\end{equation}
Furthermore, for every $\varepsilon>0$ and every fixed positive
	constant $C_\alpha$ depending only on $\alpha$, we have
	\begin{equation}\label{eq:all-mass-absorption}
		C_\alpha\int_{\mathbb R}
		\bigl(n_0^3+n_0^{3+\alpha}\bigr)dy
		\le
		\varepsilon\|\partial_yn_0\|_{L^2(\mathbb R)}^2
		+C(\alpha,M,\varepsilon).
	\end{equation}
\end{lemma}
\begin{proof}
	If $M=0$, nonnegativity and mass conservation give $n_0=0$,
	so all assertions are immediate.
	It remains to consider the case $M>0$.
	
The one-dimensional interpolation inequality gives
	\[
	\int_{\mathbb R}n_{0}^4 dy
	\le
	C\left(\int_{\mathbb R}n_{0}dy\right)^2
	\|\partial_y n_{0}\|_{L^2(\mathbb R)}^2.
	\]
	For $s\in\{2,3,3+\alpha\}$, H\"older's inequality yields
	\begin{align*}
		\int_{\mathbb R}n_{0}^sdy
		&\le
		\left(\int_{\mathbb R}n_{0}dy\right)^{\frac{4-s}{3}}
		\left(\int_{\mathbb R}n_{0}^4 dy\right)^{\frac{s-1}{3}}\le
		C_s\left(\int_{\mathbb R}n_{0}dy\right)^{\frac{s+2}{3}}
		\|\partial_yn_{0}\|_{L^2(\mathbb R)}^{\frac{2(s-1)}{3}}.
	\end{align*}
	By Lemma \ref{lem:revised-mass}, we have
	$\int_{\mathbb R}n_0(t,y)dy=\frac{M}{2\pi}$.
	Taking $s=2$ and cubing the resulting inequality proves
	\eqref{eq:ordinary-mean-GN}.
 Taking $s=3$ and $s=3+\alpha$, we similarly obtain
	\[
	C_\alpha\int_{\mathbb R}
	\bigl(n_0^3+n_0^{3+\alpha}\bigr)dy
	\le
	C(\alpha,M)
	\left[
	\bigl(\|\partial_yn_0\|_{L^2}^2\bigr)^{\frac23}
	+
	\bigl(\|\partial_yn_0\|_{L^2}^2\bigr)^{\frac{2+\alpha}{3}}
	\right].
	\]
As $0<\alpha<1$, we have
$\frac23<1$ and $\frac{2+\alpha}{3}<1$. Hence Young's inequality gives
\[
C_\alpha\int_{\mathbb R}
\bigl(n_0^3+n_0^{3+\alpha}\bigr)dy
\le
\varepsilon\|\partial_yn_0\|_{L^2(\mathbb R)}^2
+C(\alpha,M,\varepsilon).
\]
This proves \eqref{eq:all-mass-absorption}.
\end{proof}

\subsection{Estimates for the vorticity zero mode}

\begin{lemma}\label{lem:mean-vorticity-uniform}
	On every compact classical interval $[0,T]$, there holds
	\begin{equation}\label{eq:mean-vorticity-uniform}
		\sup_{0\le t\le T}\|\omega_0(t)\|_{L^\infty(\mathbb R)}
		\le
		\|(\omega_{\mathrm{in}})_0\|_{L^\infty(\mathbb R)}
		+C\sup_{0\le t\le T}\|\omega_{\neq}(t)\|_{L^2}^2.
	\end{equation}
\end{lemma}
\begin{proof}
The zero mode of the vorticity, given by $\eqref{eq:rescaled-system}_3$, satisfies 
	\[
\partial_t\omega_0-A^{-1}\partial_y^2\omega_0
=-A^{-1}\partial_y^2(v_{1,\neq}v_{2,\neq})_0,
\]
where we used $\partial_y(v_{1,\neq}v_{2,\neq})_0=(v_{2,\neq}\omega_{\neq})_0$
by incompressibility and the definition of vorticity.
Duhamel's formula therefore gives
\begin{equation}\label{Duhamel}
\omega_0(t)
=
{\rm e}^{\frac tA\partial_y^2}(\omega_{\mathrm{in}})_0
-\frac1A\int_0^t
\partial_y^2{\rm e}^{\frac{t-s}{A}\partial_y^2}
(v_{1,\neq}v_{2,\neq})_0(s)ds.
\end{equation}
By the definition of the zero mode, H\"older's inequality
and \eqref{eq:Biot-bound}, we arrive at
\begin{equation}\label{v w}
	\begin{aligned}
		\|(v_{1,\neq}v_{2,\neq})_0\|_{L^1(\mathbb R)}
		&\le C\|v_{\neq}\|_{L^2}^2\le C\|\omega_{\neq}\|_{L^2}^2,\\
		\|\partial_y(v_{1,\neq}v_{2,\neq})_0
		\|_{L^{\frac43}(\mathbb R)}&\le
		C\|v_{2,\neq}\|_{L^4}
		\|\omega_{\neq}\|_{L^2}
		\le C\|\omega_{\neq}\|_{L^2}^2.
	\end{aligned}
\end{equation}
The one-dimensional heat kernel is
$\frac1{\sqrt{4\pi\tau}}
\exp\left(-\frac{y^2}{4\tau}\right)$ for $\tau>0$.
Scaling gives
\[
\left\|
\partial_y\left[
\frac1{\sqrt{4\pi\tau}}
\exp\left(-\frac{y^2}{4\tau}\right)
\right]
\right\|_{L^4(\mathbb R)}
\le C\tau^{-\frac78},\quad
\left\|
\partial_y^2\left[
\frac1{\sqrt{4\pi\tau}}
\exp\left(-\frac{y^2}{4\tau}\right)
\right]
\right\|_{L^\infty(\mathbb R)}
\le C\tau^{-\frac32}.
\]
Young's convolution inequality consequently yields
\begin{equation}\label{heat estimates}
	\|\partial_y{\rm e}^{\tau\partial_y^2}h
	\|_{L^\infty(\mathbb R)}
	\le
	C\tau^{-\frac78}\|h\|_{L^{\frac43}(\mathbb R)},
	\qquad
	\|\partial_y^2{\rm e}^{\tau\partial_y^2}h
	\|_{L^\infty(\mathbb R)}
	\le
	C\tau^{-\frac32}\|h\|_{L^1(\mathbb R)}.
\end{equation}
Fix $t\in(0,T]$. We distinguish the cases
$t\le A$ and $t>A$.

\noindent\textbf{Case 1: $0<t\le A$.}
For every $s\in[0,t)$, we have
$0<\frac{t-s}{A}\le\frac tA\le1$.
Thus the entire time integral can be estimated
using the derivative of the averaged product.
Since spatial differentiation commutes with the
heat semigroup, setting $\tau=\frac{t-s}{A}$ and
using $\eqref{v w}_2$ and the first estimate in
\eqref{heat estimates}, we obtain
\begin{align*}
	\left\|
	\partial_y^2{\rm e}^{\tau\partial_y^2}
	(v_{1,\neq}v_{2,\neq})_0(s)
	\right\|_{L^\infty(\mathbb R)}&=
	\left\|
	\partial_y{\rm e}^{\tau\partial_y^2}
	\partial_y(v_{1,\neq}v_{2,\neq})_0(s)
	\right\|_{L^\infty(\mathbb R)}\\&\le
	C\tau^{-\frac78}
	\|\partial_y(v_{1,\neq}v_{2,\neq})_0(s)
	\|_{L^{\frac43}(\mathbb R)}\le
	C\tau^{-\frac78}\|\omega_{\neq}(s)\|_{L^2}^2.
\end{align*}
Integrating this estimate and changing variables
from $s$ to $\tau=\frac{t-s}{A}$ give
\begin{align*}
	&\frac1A\int_0^t
	\left\|
	\partial_y^2{\rm e}^{\frac{t-s}{A}\partial_y^2}
	(v_{1,\neq}v_{2,\neq})_0(s)
	\right\|_{L^\infty(\mathbb R)}ds\le
	C\sup_{0\le s\le t}\|\omega_{\neq}(s)\|_{L^2}^2
	\frac1A\int_0^t
	\left(\frac{t-s}{A}\right)^{-\frac78}ds\\&=
	C\sup_{0\le s\le t}\|\omega_{\neq}(s)\|_{L^2}^2
	\int_0^{\frac{t}{A}}\tau^{-\frac78}d\tau\le
	8C\sup_{0\le s\le t}\|\omega_{\neq}(s)\|_{L^2}^2,
\end{align*}
where we used $\frac{t}{A}\le1$ and
$\int_0^1\tau^{-\frac78}d\tau=8$.

\noindent\textbf{Case 2: $t>A$.}
We split the time integral into the intervals
$[0,t-A]$ and $[t-A,t]$.
On $[t-A,t)$, we have
$0<\frac{t-s}{A}\le1$.
The same pointwise estimate as in Case 1 gives
\begin{align*}
	\frac1A\int_{t-A}^t
	\left\|
	\partial_y^2{\rm e}^{\frac{t-s}{A}\partial_y^2}
	(v_{1,\neq}v_{2,\neq})_0(s)
	\right\|_{L^\infty(\mathbb R)}ds\le
	8C\sup_{0\le s\le t}\|\omega_{\neq}(s)\|_{L^2}^2.
\end{align*}
On $[0,t-A]$, we have $\frac{t-s}{A}\ge1$.
We now use the $L^1$ bound for the averaged product,
keeping both derivatives on the heat kernel.
By $\eqref{v w}_1$ and the second estimate in
\eqref{heat estimates}, we obtain
\begin{align*}
	\left\|
	\partial_y^2{\rm e}^{\frac{t-s}{A}\partial_y^2}
	(v_{1,\neq}v_{2,\neq})_0(s)
	\right\|_{L^\infty(\mathbb R)}&\le
	C\left(\frac{t-s}{A}\right)^{-\frac32}
	\|(v_{1,\neq}v_{2,\neq})_0(s)\|_{L^1(\mathbb R)}\\&\le
	C\left(\frac{t-s}{A}\right)^{-\frac32}
	\|\omega_{\neq}(s)\|_{L^2}^2.
\end{align*}
Consequently,
\begin{align*}
	&\frac1A\int_0^{t-A}
	\left\|
	\partial_y^2{\rm e}^{\frac{t-s}{A}\partial_y^2}
	(v_{1,\neq}v_{2,\neq})_0(s)
	\right\|_{L^\infty(\mathbb R)}ds\le
	C\sup_{0\le s\le t}\|\omega_{\neq}(s)\|_{L^2}^2
	\frac1A\int_0^{t-A}
	\left(\frac{t-s}{A}\right)^{-\frac32}ds\\&=
	C\sup_{0\le s\le t}\|\omega_{\neq}(s)\|_{L^2}^2
	\int_1^{\frac{t}{A}}\tau^{-\frac32}d\tau\le
	2C\sup_{0\le s\le t}\|\omega_{\neq}(s)\|_{L^2}^2,
\end{align*}
where we used $\int_1^\infty\tau^{-\frac32}d\tau=2$.

Combining the two cases, we conclude that
\[
\frac1A\int_0^t
\left\|
\partial_y^2{\rm e}^{\frac{t-s}{A}\partial_y^2}
(v_{1,\neq}v_{2,\neq})_0(s)
\right\|_{L^\infty(\mathbb R)}ds
\le
C\sup_{0\le s\le t}\|\omega_{\neq}(s)\|_{L^2}^2.
\]
This together with \eqref{Duhamel} yields
\[
\|\omega_0(t)\|_{L^\infty(\mathbb R)}
\le
\|(\omega_{\mathrm{in}})_0\|_{L^\infty(\mathbb R)}
+C\sup_{0\le s\le t}\|\omega_{\neq}(s)\|_{L^2}^2,
\qquad 0\le t\le T.
\]
Taking the supremum over $t\in[0,T]$ proves
\eqref{eq:mean-vorticity-uniform}.
\end{proof}
\section{Joint density--vorticity energy estimates}
\label{sec:full-energy}
We estimate the energy $E$ defined in \eqref{eq:full-energy}.
Its equivalence with the joint $L^2$ norm was established in
\eqref{eq:full-energy-equivalence}. All calculations are made
on a compact interval $[0,T]$ of a classical solution of
\eqref{eq:rescaled-system}. The constants used to choose $A$
will be independent of $T$.

\subsection{Fourier multiplier estimates}
We establish the multiplier bounds and the dissipation
estimate used in the joint energy argument.
The multiplier in \eqref{eq:full-energy} has symbol
\begin{equation}\label{eq:bounded-symbol}
	\operatorname{sgn}(k)
	\arctan\bigl(A^{-\frac13}|k|^{-\frac13}\eta\bigr),
	\qquad k\ne0,
\end{equation}
and is zero at $k=0$.
Since the symbol is real and bounded in absolute value
by $\frac{\pi}{2}$, the multiplier is self-adjoint and has
$L^2$ operator norm at most $\frac{\pi}{2}$.
It also commutes with spatial derivatives.
For $\partial_tf+y\partial_xf=0$, we have
$\widehat f(t,k,\eta)=\widehat f(0,k,\eta+kt)$.
Changing variables in the quadratic correction replaces
$\eta$ in its symbol by $\eta-kt$.
The derivative of the shifted symbol is
\[
\partial_t\left[
\operatorname{sgn}(k)
\arctan\bigl(
A^{-\frac13}|k|^{-\frac13}(\eta-kt)
\bigr)
\right]
=
-\frac{A^{-\frac13}|k|^{\frac23}}
{1+A^{-\frac23}|k|^{-\frac23}(\eta-kt)^2}.
\]
Its absolute value is bounded by
$A^{-\frac13}|k|^{\frac23}$.
For $f(0)\in H^1$, the estimate
\[
\sum_{k\ne0}|k|^{\frac23}
\int_{\mathbb R}|\widehat f(0,k,\eta)|^2 d\eta
\le\|f(0)\|_{H^1}^2
\]
justifies differentiation under the sum and integral.
Changing variables back yields
\begin{equation}\label{eq:correction-shear}
	\begin{aligned}
		&\frac d{dt}
		\sum_{k\ne0}\int_{\mathbb R}
		\operatorname{sgn}(k)
		\arctan\bigl(
		A^{-\frac13}|k|^{-\frac13}\eta
		\bigr)
		|\widehat f(t,k,\eta)|^2 d\eta\\&=
		-\sum_{k\ne0}\int_{\mathbb R}
		\frac{A^{-\frac13}|k|^{\frac23}}
		{1+A^{-\frac23}|k|^{-\frac23}\eta^2}
		|\widehat f(t,k,\eta)|^2 d\eta.
	\end{aligned}
\end{equation}
For the full density and vorticity equations,
the calculation is justified by the coordinate change
$f(t,x+ty,y)$.
Its time derivative is
\[
\partial_t[f(t,x+ty,y)]
=
(\partial_tf+y\partial_xf)(t,x+ty,y),
\]
which belongs to $L^2$ on each compact classical
interval for $f=n$ and $f=\omega_{\neq}$.
The $L^2$ differentiation rule therefore gives
the Couette contribution in \eqref{eq:correction-shear},
together with the pairings with the remaining terms.

We next combine this contribution with ordinary diffusion.
For $0<a\le1$, $k\ne0$ and $\eta\in\mathbb R$, there holds
\begin{equation}\label{eq:fourier-coercivity}
	\frac{k^2+\eta^2}{4A}
	+\frac{aA^{-\frac13}|k|^{\frac23}}
	{1+A^{-\frac23}|k|^{-\frac23}\eta^2}
	\ge
	\frac a4 A^{-\frac13}|k|^{\frac23}.
\end{equation}
Indeed, if $|\eta|\le A^{\frac13}|k|^{\frac13}$,
the second term is at least
$\frac a2A^{-\frac13}|k|^{\frac23}$.
Otherwise,
\[
\frac{k^2+\eta^2}{4A}
\ge\frac{\eta^2}{4A}
\ge\frac14A^{-\frac13}|k|^{\frac23}
\ge\frac a4A^{-\frac13}|k|^{\frac23}.
\]
Thus \eqref{eq:fourier-coercivity} holds in both cases.

Multiplying \eqref{eq:fourier-coercivity} by
$|\widehat f(k,\eta)|^2$, integrating over $\mathbb R$
and summing over $k\ne0$, Parseval's identity gives
\[
\begin{aligned}
	&\frac1{4A}\|\nabla f_{\neq}\|_{L^2}^2
	+a\sum_{k\ne0}\int_{\mathbb R}
	\frac{A^{-\frac13}|k|^{\frac23}}
	{1+A^{-\frac23}|k|^{-\frac23}\eta^2}
	|\widehat f(k,\eta)|^2 d\eta\\&\ge
	\frac a4 A^{-\frac13}
	\sum_{k\ne0}\int_{\mathbb R}
	|k|^{\frac23}|\widehat f(k,\eta)|^2 d\eta\ge
	\frac a4 A^{-\frac13}\|f_{\neq}\|_{L^2}^2,
\end{aligned}
\]
where the last inequality uses $|k|\ge1$.
In the joint energy estimate, this inequality is applied
to the density and vorticity with
\[
a=
\frac1{
	8\pi(1+2K_\infty)^{\alpha-1}
	[1+2(1-\alpha)K_\infty]}\quad\text{and}
\quad
a=\frac1{8\pi},
\]
respectively.
Both coefficients belong to $(0,1]$.
Once the bootstrap constants are fixed,
the resulting $A^{-\frac13}$ dissipation absorbs
the non-zero-mode errors with coefficient $A^{-1}$
for sufficiently large $A$.
\begin{lemma}\label{lem:mean-transport-correction}
Let $A\ge1$ and $f=f_{\neq}\in H^1(\mathbb T\times\mathbb R)$.
Suppose that $v_{1,0}\in W^{1,\infty}(\mathbb R)$
	is real-valued. Then
	\begin{equation}\label{eq:mean-transport-correction}
		\begin{aligned}
			&\left|
			-\frac2A\operatorname{Re}
			\sum_{k\ne0}\int_{\mathbb R}
			\operatorname{sgn}(k)
			\arctan\bigl(
			A^{-\frac13}|k|^{-\frac13}\eta
			\bigr)
			\overline{\widehat f(k,\eta)}
			[v_{1,0}\partial_xf]^{\wedge}(k,\eta)
			d\eta
			\right|\\
			&\quad\le
			\frac CA
			\|\partial_yv_{1,0}\|_{L^\infty(\mathbb R)}
			\sum_{k\ne0}A^{-\frac13}|k|^{\frac23}
			\int_{\mathbb R}
			|\widehat f(k,\eta)|^2 d\eta.
		\end{aligned}
	\end{equation}
\end{lemma}
\begin{proof}
	For $k\ne0$, we use the unitary Fourier transforms
	\[
	f_k(y)=\frac1{\sqrt{2\pi}}
	\int_{\mathbb T}f(x,y){\rm e}^{-\mathrm ikx}dx,
	\qquad
	\widehat f(k,\eta)=\frac1{\sqrt{2\pi}}
	\int_{\mathbb R}f_k(y){\rm e}^{-\mathrm i\eta y}dy.
	\]
The analysis is based on the integral representation
	\begin{equation}\label{eq:mean-transport-arctan}
		\arctan\bigl(A^{-\frac13}|k|^{-\frac13}\eta\bigr)
		=
		\int_0^\infty
		{\rm e}^{-A^{\frac13}|k|^{\frac13}r}
		\frac{\sin(\eta r)}r dr.
	\end{equation}
	Indeed, both sides vanish at $\eta=0$, and
	differentiating the right-hand side gives
	\[
	\int_0^\infty
	{\rm e}^{-A^{\frac13}|k|^{\frac13}r}
	\cos(\eta r)dr
	=\frac{A^{\frac13}|k|^{\frac13}}
	{A^{\frac23}|k|^{\frac23}+\eta^2},
	\]
	which is also the derivative of the left-hand side.
	Using \eqref{eq:mean-transport-arctan} and Fourier
	inversion, we obtain
	\begin{equation}\label{eq:mean-transport-kernel}
		\begin{aligned}
			&\frac1{\sqrt{2\pi}}\int_{\mathbb R}
			{\rm e}^{\mathrm i\eta y}
			\operatorname{sgn}(k)
			\arctan\bigl(A^{-\frac13}|k|^{-\frac13}\eta\bigr)
			\widehat f(k,\eta)d\eta\\=~&
			\frac{\operatorname{sgn}(k)}{2\mathrm i}
			\int_0^\infty
			{\rm e}^{-A^{\frac13}|k|^{\frac13}r}
			\frac{f_k(y+r)-f_k(y-r)}rdr\\=~&
			\frac{\mathrm i\operatorname{sgn}(k)}2
			\operatorname{p.v.}\int_{\mathbb R}
			\frac{{\rm e}^{-A^{\frac13}|k|^{\frac13}|y-z|}
			}{y-z}f_k(z)dz.
		\end{aligned}
	\end{equation}
	To estimate the expression on the left-hand side of
\eqref{eq:mean-transport-correction}, consider
$\partial_tf=-A^{-1}v_{1,0}\partial_xf$.
Then
$\partial_t f_k=-\frac{\mathrm ik}{A}v_{1,0}f_k$,
and the expression inside the absolute value is
the time derivative of the quadratic correction. 
	Substituting it when
	differentiating the correction, and using Parseval's
	identity and the kernel representation
	\eqref{eq:mean-transport-kernel}, gives
	\begin{equation}\label{eq:mean-transport-pairing}
		\begin{aligned}
			&\frac d{dt}\int_{\mathbb R}
			\operatorname{sgn}(k)
			\arctan\bigl(A^{-\frac13}|k|^{-\frac13}\eta\bigr)
			|\widehat f(k,\eta)|^2 d\eta\\=~&
			2\operatorname{Re}\int_{\mathbb R}
			\operatorname{sgn}(k)
			\arctan\bigl(A^{-\frac13}|k|^{-\frac13}\eta\bigr)
			\overline{\widehat f(k,\eta)}
			\partial_t\widehat f(k,\eta)d\eta\\=~&
			\frac{|k|}{A}\operatorname{Re}
			\lim_{\varepsilon\to 0}
			\int_{\mathbb R}
			\int_{\{z\in\mathbb R:\,|y-z|>\varepsilon\}}
			\frac{
				{\rm e}^{-A^{\frac13}|k|^{\frac13}|y-z|}
			}{y-z}
			v_{1,0}(z)\overline{f_k(y)}f_k(z)dzdy.
		\end{aligned}
	\end{equation}
For each $\varepsilon>0$, interchange $y$ and $z$
	in the truncated double integral in
	\eqref{eq:mean-transport-pairing}.
	Since $v_{1,0}$ is real,
	$\operatorname{Re}(\overline{f_k(y)}f_k(z))$
	is unchanged under this interchange, whereas
	$(y-z)^{-1}$ changes sign.
	Averaging the original and interchanged expressions
	therefore gives
	\begin{equation}\label{eq:mean-transport-symmetrized}
		\begin{aligned}
			&\frac d{dt}\int_{\mathbb R}
			\operatorname{sgn}(k)
			\arctan\bigl(A^{-\frac13}|k|^{-\frac13}\eta\bigr)
			|\widehat f(k,\eta)|^2 d\eta\\
			&\quad=
			-\frac{|k|}{2A}\operatorname{Re}
			\int_{\mathbb R}\int_{\mathbb R}
			{\rm e}^{-A^{\frac13}|k|^{\frac13}|y-z|}
			\frac{v_{1,0}(y)-v_{1,0}(z)}{y-z}
			\overline{f_k(y)}f_k(z)dzdy.
		\end{aligned}
	\end{equation}
	The passage to the limit is justified by the
	Lipschitz bound
	\[
	\left|
	\frac{v_{1,0}(y)-v_{1,0}(z)}{y-z}
	\right|
	\le
	\|\partial_yv_{1,0}\|_{L^\infty(\mathbb R)}\quad {\text{for}}\quad
	y\ne z
	\]
	and the following integrability estimate:
	\begin{equation}\label{eq:mean-transport-kernel-bound}
		\begin{aligned}
			&\int_{\mathbb R}\int_{\mathbb R}
			{\rm e}^{-A^{\frac13}|k|^{\frac13}|y-z|}
			|f_k(y)||f_k(z)|dzdy\\\le~&
			\frac12\int_{\mathbb R}\int_{\mathbb R}
			{\rm e}^{-A^{\frac13}|k|^{\frac13}|y-z|}
			\bigl(|f_k(y)|^2+|f_k(z)|^2\bigr)dzdy=
			2A^{-\frac13}|k|^{-\frac13}
			\|f_k\|_{L^2(\mathbb R)}^2.
		\end{aligned}
	\end{equation}
	In the last equality, we used
	$\int_{\mathbb R}
	{\rm e}^{-A^{\frac13}|k|^{\frac13}|y-z|}dz
	=2A^{-\frac13}|k|^{-\frac13}$.
	Thus the double integral in
	\eqref{eq:mean-transport-symmetrized}
	is absolutely convergent.
	Combining \eqref{eq:mean-transport-symmetrized}
	with \eqref{eq:mean-transport-kernel-bound}, we obtain
		\begin{align*}
			\left|
			\frac d{dt}\int_{\mathbb R}
			\operatorname{sgn}(k)
			\arctan\bigl(A^{-\frac13}|k|^{-\frac13}\eta\bigr)
			|\widehat f(k,\eta)|^2 d\eta
			\right|\le
			\frac1A
			\|\partial_yv_{1,0}\|_{L^\infty(\mathbb R)}
			A^{-\frac13}|k|^{\frac23}
			\|f_k\|_{L^2(\mathbb R)}^2.
		\end{align*}
	By Parseval's identity, $A\ge1$ and $|k|\ge1$, there holds
\[
\begin{aligned}
	\sum_{k\ne0}A^{-\frac13}|k|^{\frac23}
	\|f_k\|_{L^2(\mathbb R)}^2\le
	\sum_{k\ne0}k^2\|f_k\|_{L^2(\mathbb R)}^2=
	\|\partial_xf\|_{L^2(\mathbb T\times\mathbb R)}^2
	<\infty.
\end{aligned}
\]
Thus the sum of the preceding mode estimates
is absolutely convergent.
Summing over $k\ne0$ and applying Parseval's identity
in $y$ proves \eqref{eq:mean-transport-correction}.
\end{proof}

\subsection{The joint energy estimate}
We now estimate the energy defined in
\eqref{eq:full-energy} under the bootstrap bounds
\eqref{eq:full-bootstrap}. The constants $K_{\neq}$ and
$K_\infty$ remain fixed throughout the estimates.
\begin{proposition}\label{prop:joint-energy}
	Under \eqref{eq:full-bootstrap}, there exists a threshold
	for $A$, depending only on $\alpha$, $M$, $K_{\neq}$,
	$K_\infty$ and
	$\|\omega_{\mathrm{in},0}\|_{L^\infty(\mathbb R)}$,
	such that
	\begin{equation}
			\begin{aligned}\label{eq:joint-differential}
			&E'(t)+
			\frac1{
				64\pi(1+2K_\infty)^{\alpha-1}
				[1+2(1-\alpha)K_\infty]}
			\sum_{k\ne0}\int_{\mathbb R}
			A^{-\frac13}|k|^{\frac23}
			|\widehat n(t,k,\eta)|^2 d\eta\\&+\frac1A
			\left(
			\|\nabla n\|_{L^2}^2
			+\|\nabla\omega_{\neq}\|_{L^2}^2
			\right)+\frac1{64\pi}
			\sum_{k\ne0}\int_{\mathbb R}
			A^{-\frac13}|k|^{\frac23}
			|\widehat\omega(t,k,\eta)|^2 d\eta
			\le\frac{C(\alpha,M)}A
		\end{aligned}
	\end{equation}
	for almost every $t\in[0,T]$.
	The constant on the right is independent of
	$A$, $T$, $K_{\neq}$ and $K_\infty$.
\end{proposition}
\begin{proof}
	The bootstrap bounds \eqref{eq:full-bootstrap},
	nonnegativity and \eqref{eq:full-energy-equivalence}
	imply
	\begin{equation}\label{eq:joint-proof-bootstrap-consequences}
		\begin{aligned}
			\|n(t)\|_{L^2}^2+\|\omega_{\neq}(t)\|_{L^2}^2\le\frac83K_{\neq},\quad
			0\le n_0,n\le2K_\infty,
			\quad |n_{\neq}|\le2K_\infty.
		\end{aligned}
	\end{equation}
	By Lemma \ref{lem:mean-vorticity-uniform}, there holds
	\begin{equation}\label{eq:joint-proof-mean-bound}
		\|\omega_0\|_{L^\infty(0,T;L^\infty(\mathbb R))}
		\le
		\|(\omega_{\mathrm{in}})_0\|_{L^\infty(\mathbb R)}
		+CK_{\neq}.
	\end{equation}
	Moreover, Lemma \ref{lem:revised-coefficients} gives
	\begin{equation}\label{eq:joint-diffusion-range}
		1\le (1+n)^{\alpha-1}[1+(1-\alpha)n]
		\le (1+2K_\infty)^{\alpha-1}
		[1+2(1-\alpha)K_\infty].
	\end{equation}
	We also record a bound used in the Fourier estimates.
	Since the symbol in \eqref{eq:bounded-symbol}
	has absolute value at most $\frac{\pi}{2}$,
	Parseval's identity gives, for $f=n$ or
	$f=\omega_{\neq}$,
	\begin{equation}\label{eq:joint-multiplier-gradient}
		\begin{aligned}
			&\sum_{k\ne0}\int_{\mathbb R}
			(k^2+\eta^2)
			\left|
			\arctan\bigl(A^{-\frac13}|k|^{-\frac13}\eta\bigr)
			\right|^2
			|\widehat f(k,\eta)|^2 d\eta\\&\le
			\frac{\pi^2}{4}
			\sum_{k\ne0}\int_{\mathbb R}
			(k^2+\eta^2)|\widehat f(k,\eta)|^2 d\eta
			=\frac{\pi^2}{4}\|\nabla f_{\neq}\|_{L^2}^2.
		\end{aligned}
	\end{equation}

\noindent{\bf Step I. The density energy.}
	{\bf (i) The ordinary $L^{2}$ estimate.}
	Testing $\eqref{eq:rescaled-system}_1$ with $2n$, and using integration by parts together with $\nabla\cdot v=0$,
	gives
	\begin{equation}\label{eq:joint-density-L2-identity}
		\begin{aligned}
			&\frac d{dt}\|n\|_{L^2}^2
			+\frac2A\int_{\mathbb T\times\mathbb R}
			(1+n)^{\alpha-1}[1+(1-\alpha)n]
			|\nabla n|^2 dxdy\\&=
			\frac2A\int_{\mathbb T\times\mathbb R}
			n(1+n)^\alpha\nabla c\cdot\nabla n dxdy.
		\end{aligned}
	\end{equation}
	Due to
	$\nabla\left(
	\int_0^{n(t,x,y)}s(1+s)^\alpha ds\right)=
	n(1+n)^\alpha\nabla n$,
	integration by parts and $\Delta c=c-n$ yield
	\begin{equation}\label{eq:joint-chemotactic-production}
		\begin{aligned}
			&\frac2A\int_{\mathbb T\times\mathbb R}
			n(1+n)^\alpha\nabla c\cdot\nabla n dxdy=
			-\frac2A\int_{\mathbb T\times\mathbb R}
			\left(
			\int_0^{n(t,x,y)}s(1+s)^\alpha ds
			\right)\Delta c dxdy\\
			&=\frac2A\int_{\mathbb T\times\mathbb R}
			\left(
			\int_0^{n(t,x,y)}s(1+s)^\alpha ds
			\right)(n-c)dxdy\le
			\frac{C_\alpha}{A}
			\int_{\mathbb T\times\mathbb R}
			\bigl(n^3+n^{3+\alpha}\bigr)dxdy.
		\end{aligned}
	\end{equation}
	The last inequality follows from $c\ge0$ and
	\[
	0\le
	\int_0^{n(t,x,y)}s(1+s)^\alpha ds
	\le
	\int_0^{n(t,x,y)}(s+s^{1+\alpha}) ds
	=
	\frac{n(t,x,y)^2}{2}
	+\frac{n(t,x,y)^{2+\alpha}}{2+\alpha}.
	\]
For $s\in\{3,3+\alpha\}$, the decomposition
$n=n_0+n_{\neq}$ and
\eqref{eq:joint-proof-bootstrap-consequences} give
\[
\begin{aligned}
	\int_{\mathbb T\times\mathbb R}n^s dxdy&\le
	2^{s-1}\left(
	2\pi\int_{\mathbb R}n_0^s dy
	+\int_{\mathbb T\times\mathbb R}|n_{\neq}|^s dxdy
	\right)\\&\le
	C_s\int_{\mathbb R}n_0^s dy
	+C_s(2K_\infty)^{s-2}\|n_{\neq}\|_{L^2}^2.
\end{aligned}
\]
Summing these estimates and applying
\eqref{eq:all-mass-absorption} to the resulting
zero-mode terms with $\varepsilon=\frac{\pi}{16}$,
we obtain
\[
\begin{aligned}
	\frac{C_\alpha}{A}
	\int_{\mathbb T\times\mathbb R}
	\bigl(n^3+n^{3+\alpha}\bigr)dxdy\le~&
	\frac{\pi}{16A}
	\|\partial_yn_0\|_{L^2(\mathbb R)}^2
	+\frac{C(\alpha,M)}A
	+\frac{C(\alpha,K_\infty)}A
	\|n_{\neq}\|_{L^2}^2\\\le~&
	\frac1{32A}\|\nabla n\|_{L^2}^2
	+\frac{C(\alpha,M)}A
	+\frac{C(\alpha,K_\infty)}A
	\|n_{\neq}\|_{L^2}^2.
\end{aligned}
\]
The last inequality follows from
$\|\nabla n\|_{L^2}^2=
2\pi\|\partial_yn_0\|_{L^2(\mathbb R)}^2
+\|\nabla n_{\neq}\|_{L^2}^2$.
Combining this bound with \eqref{eq:joint-diffusion-range}, \eqref{eq:joint-density-L2-identity} and
\eqref{eq:joint-chemotactic-production}
 yields
	\begin{equation}\label{eq:joint-density-basic-bound}
		\frac d{dt}\|n\|_{L^2}^2
		+\frac{63}{32A}\|\nabla n\|_{L^2}^2
		\le
		\frac{C(\alpha,M)}A
		+\frac{C(\alpha,K_\infty)}A
		\|n_{\neq}\|_{L^2}^2.
	\end{equation}

	\noindent{\bf (ii) The Fourier correction estimate.}
We compute the time derivative of the density
Fourier correction in \eqref{eq:full-energy}.
Taking the Fourier transform of
$\eqref{eq:rescaled-system}_1$, for $k\ne0$, we obtain
\begin{equation}\label{eq:joint-density-fourier-equation}
	\begin{aligned}
		\partial_t\widehat n(t,k,\eta)
		-k\partial_\eta\widehat n(t,k,\eta)
		={}&\frac1A
		\left[
		\nabla\cdot\left(
		(1+n)^{\alpha-1}[1+(1-\alpha)n]\nabla n
		\right)
		\right]^{\wedge}(t,k,\eta)\\
		&-\frac1A
		\left[
		\nabla\cdot\bigl(n(1+n)^\alpha\nabla c+vn\bigr)
		\right]^{\wedge}(t,k,\eta).
	\end{aligned}
\end{equation}
We multiply \eqref{eq:joint-density-fourier-equation}
by $2\overline{\widehat n(t,k,\eta)}$
and take the real part.
For the time derivative and Couette terms, there holds
\[
\begin{aligned}
	2\operatorname{Re}\left(
	\overline{\widehat n}\,\partial_t\widehat n
	\right)=\partial_t|\widehat n|^2,\quad
	2\operatorname{Re}\left(
	\overline{\widehat n}\,k\partial_\eta\widehat n
	\right)=k\partial_\eta|\widehat n|^2.
\end{aligned}
\]
For the terms containing spatial derivatives,
we use $\overline{[\nabla n]^{\wedge}(t,k,\eta)}
=-\mathrm i(k,\eta)\,
\overline{\widehat n(t,k,\eta)}$ and
	\begin{align*}
		&\left[
		\nabla\cdot\left(
		(1+n)^{\alpha-1}[1+(1-\alpha)n]\nabla n
		\right)
		\right]^{\wedge}(t,k,\eta)\\
		&\qquad=
		\mathrm i(k,\eta)\cdot
		\left[
		(1+n)^{\alpha-1}[1+(1-\alpha)n]\nabla n
		\right]^{\wedge}(t,k,\eta).
	\end{align*}
Thus the diffusion term satisfies
\[
\begin{aligned}
	&\overline{\widehat n(t,k,\eta)}
	\left[
	\nabla\cdot\left(
	(1+n)^{\alpha-1}[1+(1-\alpha)n]\nabla n
	\right)
	\right]^{\wedge}(t,k,\eta)\\
	&\qquad=
	-\overline{[\nabla n]^{\wedge}(t,k,\eta)}
	\cdot
	\left[
	(1+n)^{\alpha-1}[1+(1-\alpha)n]\nabla n
	\right]^{\wedge}(t,k,\eta).
\end{aligned}
\]
Applying the same differentiation identities to
the chemotactic and velocity terms gives the
opposite sign because their divergence appears
with a minus sign in
\eqref{eq:joint-density-fourier-equation}.
Moreover, at $k\ne0$, inserting $P_{\neq}$
does not change the Fourier coefficient.
Consequently,
\begin{equation}\label{eq:joint-density-fourier-square}
	\begin{aligned}
		&\partial_t|\widehat n(t,k,\eta)|^2
		-k\partial_\eta|\widehat n(t,k,\eta)|^2\\=~&
		-\frac2A\operatorname{Re}\left\{
		\overline{[\nabla n]^{\wedge}(t,k,\eta)}
		\cdot
		\left[
		(1+n)^{\alpha-1}[1+(1-\alpha)n]\nabla n
		\right]^{\wedge}(t,k,\eta)
		\right\}\\
		&+
		\frac2A\operatorname{Re}\left\{
		\overline{[\nabla n]^{\wedge}(t,k,\eta)}
		\cdot
		\left[
		P_{\neq}\bigl(n(1+n)^\alpha\nabla c+vn\bigr)
		\right]^{\wedge}(t,k,\eta)
		\right\}.
	\end{aligned}
\end{equation}
We now multiply \eqref{eq:joint-density-fourier-square}
by
$\operatorname{sgn}(k)
\arctan\bigl(A^{-\frac13}|k|^{-\frac13}\eta\bigr)$,
integrate over $\eta\in\mathbb R$,
and sum over $k\ne0$.
The time derivative gives the derivative of
the quadratic Fourier expression.
For the Couette term, integration by parts in $\eta$ gives
\begin{equation}\label{eq:joint-density-couette-calculation}
	\begin{aligned}
		&-\sum_{k\ne0}\int_{\mathbb R}
		\operatorname{sgn}(k)
		\arctan\bigl(A^{-\frac13}|k|^{-\frac13}\eta\bigr)
		k\partial_\eta|\widehat n(t,k,\eta)|^2 d\eta\\
		&=
		\sum_{k\ne0}\int_{\mathbb R}
		k\partial_\eta\left[
		\operatorname{sgn}(k)
		\arctan\bigl(A^{-\frac13}|k|^{-\frac13}\eta\bigr)
		\right]
		|\widehat n(t,k,\eta)|^2 d\eta\\&=
		\sum_{k\ne0}\int_{\mathbb R}
		\frac{
			A^{-\frac13}|k|^{\frac23}
			|\widehat n(t,k,\eta)|^2
		}{
			1+A^{-\frac23}|k|^{-\frac23}\eta^2
		}d\eta.
	\end{aligned}
\end{equation}
Finally, multiply the resulting identity by
the density correction coefficient
$$
\frac1{
	8\pi(1+2K_\infty)^{\alpha-1}
	[1+2(1-\alpha)K_\infty]}.$$
Using \eqref{eq:joint-density-fourier-square}
and \eqref{eq:joint-density-couette-calculation},
we obtain
\begin{align}
	&\frac1{
		8\pi(1+2K_\infty)^{\alpha-1}
		[1+2(1-\alpha)K_\infty]}
	\frac d{dt}
	\sum_{k\ne0}\int_{\mathbb R}
	\operatorname{sgn}(k)
	\arctan\bigl(A^{-\frac13}|k|^{-\frac13}\eta\bigr)
	|\widehat n(t,k,\eta)|^2 d\eta
	\notag\\&+
	\frac1{
		8\pi(1+2K_\infty)^{\alpha-1}
		[1+2(1-\alpha)K_\infty]}
	\sum_{k\ne0}\int_{\mathbb R}
	\frac{
		A^{-\frac13}|k|^{\frac23}
		|\widehat n(t,k,\eta)|^2
	}{
		1+A^{-\frac23}|k|^{-\frac23}\eta^2
	}\,d\eta
	\notag\\&=
	-\frac2{
		8\pi A(1+2K_\infty)^{\alpha-1}
		[1+2(1-\alpha)K_\infty]}
	\operatorname{Re}
	\sum_{k\ne0}\int_{\mathbb R}
	\operatorname{sgn}(k)
	\arctan\bigl(A^{-\frac13}|k|^{-\frac13}\eta\bigr)
	\notag\\&\times
	\overline{[\nabla n]^{\wedge}(t,k,\eta)}
	\cdot
	\left[
	(1+n)^{\alpha-1}[1+(1-\alpha)n]\nabla n
	\right]^{\wedge}(t,k,\eta)d\eta
	\notag\\&+
	\frac2{
		8\pi A(1+2K_\infty)^{\alpha-1}
		[1+2(1-\alpha)K_\infty]}
	\operatorname{Re}
	\sum_{k\ne0}\int_{\mathbb R}
	\operatorname{sgn}(k)
	\arctan\bigl(A^{-\frac13}|k|^{-\frac13}\eta\bigr)
	\notag\\&\times
	\overline{[\nabla n]^{\wedge}(t,k,\eta)}
	\cdot\left[
	P_{\neq}\left(n(1+n)^\alpha\nabla c+vn\right)
	\right]^{\wedge}(t,k,\eta)d\eta.
	\label{eq:joint-density-correction-identity}
\end{align}
We next estimate the diffusion, chemotactic
and velocity terms on the right-hand side
of \eqref{eq:joint-density-correction-identity}
separately. 
For the diffusion term, Cauchy--Schwarz, Parseval's identity, 	\eqref{eq:joint-diffusion-range} and
	\eqref{eq:joint-multiplier-gradient}
	bound its absolute value by
	\begin{equation}\label{eq:joint-density-diffusion-correction}
		\begin{aligned}
			\frac1{8A}
			\|\nabla n_{\neq}\|_{L^2}\|\nabla n\|_{L^2}
			\le
			\frac1{8A}\|\nabla n\|_{L^2}^2.
		\end{aligned}
	\end{equation}
	For the chemotactic term, we first estimate
	$P_{\neq}(n(1+n)^\alpha\nabla c)$.
	Since $n_0(1+n_0)^\alpha\nabla c_0$ is independent
	of $x$, we have
	\begin{equation}\label{eq:joint-chemotactic-flux-decomposition}
		\begin{aligned}
			P_{\neq}\bigl(n(1+n)^\alpha\nabla c\bigr)
			=P_{\neq}\Bigl([n(1+n)^\alpha-n_0(1+n_0)^\alpha]\nabla c_0+n(1+n)^\alpha\nabla c_{\neq}
			\Bigr).
		\end{aligned}
	\end{equation}
	On $[0,2K_\infty]$, the direct calculation implies
	\[
	\frac d{ds}[s(1+s)^\alpha]
	=
	(1+s)^\alpha+\alpha s(1+s)^{\alpha-1}
	\le
	(1+\alpha)(1+2K_\infty)^\alpha.
	\]
	Thus the mean value theorem gives
	\[
	|n(1+n)^\alpha-n_0(1+n_0)^\alpha|
	\le
	(1+\alpha)(1+2K_\infty)^\alpha|n_{\neq}|.
	\]
	Using this bound,
	$n(1+n)^\alpha\le2K_\infty(1+2K_\infty)^\alpha$, \eqref{eq:nonzero-chemical-gradient},
	\eqref{eq:mass-chemical} and
	\eqref{eq:joint-chemotactic-flux-decomposition},
	we obtain
	\begin{equation}\label{eq:joint-chemotactic-flux-bound}
		\begin{aligned}
			\bigl\|
			P_{\neq}\bigl(n(1+n)^\alpha\nabla c\bigr)&
			\bigr\|_{L^2}\le
			(1+\alpha)(1+2K_\infty)^\alpha
			\frac M{4\pi}\|n_{\neq}\|_{L^2}
			\\&+
			2K_\infty(1+2K_\infty)^\alpha
			\|\nabla c_{\neq}\|_{L^2}\le
			C(\alpha,M,K_\infty)\|n_{\neq}\|_{L^2}.
		\end{aligned}
	\end{equation}
	Apply Cauchy--Schwarz and
	\eqref{eq:joint-multiplier-gradient}
	to the part of \eqref{eq:joint-density-correction-identity}
	containing this product.
	By \eqref{eq:joint-chemotactic-flux-bound},
	its absolute value is bounded by
	\begin{equation}\label{eq:joint-corrected-chemotaxis-bound}
		\begin{aligned}
			&\frac{
				\|\nabla n_{\neq}\|_{L^2}
				\|P_{\neq}(n(1+n)^\alpha\nabla c)\|_{L^2}
			}{
				8A(1+2K_\infty)^{\alpha-1}
				[1+2(1-\alpha)K_\infty]
			}\\\le~&
			\frac{C(\alpha,M,K_\infty)}A
			\|n_{\neq}\|_{L^2}\|\nabla n\|_{L^2}\le
			\frac1{32A}\|\nabla n\|_{L^2}^2
			+\frac{C(\alpha,M,K_\infty)}A
			\|n_{\neq}\|_{L^2}^2.
		\end{aligned}
	\end{equation}
	For the velocity term, $v_{2,0}=0$ gives
	$-\frac1A\nabla\cdot(vn)=
	-\frac1Av_{1,0}\partial_xn_{\neq}
	-\frac1A\nabla\cdot(v_{\neq}n)$.
	We treat these two terms separately.
	For the zero-mode velocity term, apply
	\eqref{eq:mean-transport-correction}
	with $f=n_{\neq}$ and multiply by the density
	correction coefficient.
	Using $\partial_yv_{1,0}=\omega_0$ and
	\eqref{eq:joint-proof-mean-bound},
	we obtain the absolute bound
	\begin{equation}\label{eq:joint-density-mean-transport-bound}
		\begin{aligned}
			\frac{C\sum_{k\ne0}\int_{\mathbb R}
				A^{-\frac13}|k|^{\frac23}
				|\widehat n(t,k,\eta)|^2 d\eta}{
				8\pi A(1+2K_\infty)^{\alpha-1}
				[1+2(1-\alpha)K_\infty]}
			\left(
			\|(\omega_{\mathrm{in}})_0\|_{L^\infty(\mathbb R)}
			+CK_{\neq}
			\right).
		\end{aligned}
	\end{equation}
	For the non-zero mode velocity term,
	\eqref{eq:Biot-bound} implies
	\[
	\|P_{\neq}(v_{\neq}n)\|_{L^2}
	\le
	\|n\|_{L^\infty}\|v_{\neq}\|_{L^2}
	\le
	CK_\infty\|\omega_{\neq}\|_{L^2}.
	\]
	Applying this bound and
	\eqref{eq:joint-multiplier-gradient}
	to its gradient pairing in
	\eqref{eq:joint-density-correction-identity}
	gives
	\begin{equation}\label{eq:joint-density-velocity-bound}
		\begin{aligned}
			&\frac{
				\|\nabla n_{\neq}\|_{L^2}
				\|P_{\neq}(v_{\neq}n)\|_{L^2}
			}{
				8A(1+2K_\infty)^{\alpha-1}
				[1+2(1-\alpha)K_\infty]
			}\\\le~&
			\frac{CK_\infty}{A}
			\|\omega_{\neq}\|_{L^2}\|\nabla n\|_{L^2}\le
			\frac1{32A}\|\nabla n\|_{L^2}^2
			+\frac{C(K_\infty)}A
			\|\omega_{\neq}\|_{L^2}^2.
		\end{aligned}
	\end{equation}
Combining the ordinary density estimate \eqref{eq:joint-density-basic-bound} with the correction identity \eqref{eq:joint-density-correction-identity}, and using the estimates \eqref{eq:joint-density-diffusion-correction}, \eqref{eq:joint-corrected-chemotaxis-bound}, \eqref{eq:joint-density-mean-transport-bound}, and \eqref{eq:joint-density-velocity-bound}, we obtain
	\begin{align}
		&\frac d{dt}\left[
		\|n\|_{L^2}^2
		+\frac{\sum_{k\ne0}\int_{\mathbb R}
			\operatorname{sgn}(k)
			\arctan\bigl(A^{-\frac13}|k|^{-\frac13}\eta\bigr)
			|\widehat n(t,k,\eta)|^2 d\eta}{
			8\pi(1+2K_\infty)^{\alpha-1}
			[1+2(1-\alpha)K_\infty]}
		\right]
		\notag\\
		&+
		\frac{57}{32A}\|\nabla n\|_{L^2}^2
		+
		\frac1{
			8\pi(1+2K_\infty)^{\alpha-1}
			[1+2(1-\alpha)K_\infty]}
		\sum_{k\ne0}\int_{\mathbb R}
		\frac{
			A^{-\frac13}|k|^{\frac23}
			|\widehat n(t,k,\eta)|^2
		}{
			1+A^{-\frac23}|k|^{-\frac23}\eta^2
		}d\eta\notag\\&\le
		\frac{C(\alpha,M)}A
		+\frac{C(\alpha,M,K_\infty)}A
		\|n_{\neq}\|_{L^2}^2
		+\frac{C(K_\infty)}A
		\|\omega_{\neq}\|_{L^2}^2
		\notag\\
		&+
		\frac{C	\left(
			\|(\omega_{\mathrm{in}})_0\|_{L^\infty(\mathbb R)}
			+CK_{\neq}
			\right)
			\sum_{k\ne0}\int_{\mathbb R}
			A^{-\frac13}|k|^{\frac23}
			|\widehat n(t,k,\eta)|^2 d\eta}{
			8\pi A(1+2K_\infty)^{\alpha-1}
			[1+2(1-\alpha)K_\infty]}.
		\label{eq:joint-density-energy-bound}
	\end{align}

	\noindent{\bf Step II. The vorticity energy.}
	Applying $P_{\neq}$ to $\eqref{eq:rescaled-system}_3$
	and using $v_{2,0}=0$, we obtain
	\begin{equation}\label{eq:joint-vorticity-equation}
		\begin{aligned}
			\partial_t\omega_{\neq}
			+y\partial_x\omega_{\neq}
			-\frac{\Delta\omega_{\neq}}A
			=-\frac1A\partial_xn_{\neq}
			-\frac1A v_{1,0}\partial_x\omega_{\neq}-\frac1A P_{\neq}\nabla\cdot
			\bigl(v_{\neq}\omega_0
			+v_{\neq}\omega_{\neq}\bigr).
		\end{aligned}
	\end{equation}
	
	\noindent{\bf (i) The ordinary $L^2$ estimate.}
	Testing \eqref{eq:joint-vorticity-equation}
	with $2\omega_{\neq}$, and using integration by parts together with $\nabla\cdot v=0$, gives
	\begin{equation}\label{eq:joint-vorticity-L2-identity}
		\begin{aligned}
			&\frac d{dt}\|\omega_{\neq}\|_{L^2}^2
			+\frac{2}{A}\|\nabla\omega_{\neq}\|_{L^2}^2\\=~&
			\frac2A\int_{\mathbb T\times\mathbb R}
			n_{\neq}\partial_x\omega_{\neq}dxdy
			+\frac2A\int_{\mathbb T\times\mathbb R}
			P_{\neq}(v_{\neq}\omega_0)
			\cdot\nabla\omega_{\neq}dxdy.
		\end{aligned}
	\end{equation}
	For the first term on the right-hand side, using
	Cauchy--Schwarz and Young's inequalities, we obtain
	\begin{equation}\label{eq:joint-vorticity-density-bound}
		\begin{aligned}
			\left|
			\frac2A\int_{\mathbb T\times\mathbb R}
			n_{\neq}\partial_x\omega_{\neq}dxdy
			\right|\le
			\frac2A\|n_{\neq}\|_{L^2}
			\|\nabla\omega_{\neq}\|_{L^2}\le
			\frac1{64A}\|\nabla\omega_{\neq}\|_{L^2}^2
			+\frac CA\|n_{\neq}\|_{L^2}^2.
		\end{aligned}
	\end{equation}
	For the second term, \eqref{eq:Biot-bound} and
	\eqref{eq:joint-proof-mean-bound} imply
		\begin{align}\label{eq:joint-vorticity-mean-product}
			\|P_{\neq}(v_{\neq}\omega_0)\|_{L^2}\le
			\|\omega_0\|_{L^\infty(\mathbb R)}
			\|v_{\neq}\|_{L^2}\le
			C\left(
			\|(\omega_{\mathrm{in}})_0\|_{L^\infty(\mathbb R)}
			+CK_{\neq}
			\right)\|\omega_{\neq}\|_{L^2}.
		\end{align}
	Thus
		\begin{align*}
			&\left|
			\frac2A\int_{\mathbb T\times\mathbb R}
			P_{\neq}(v_{\neq}\omega_0)
			\cdot\nabla\omega_{\neq}dxdy
			\right|\le
			\frac2A
			\|P_{\neq}(v_{\neq}\omega_0)\|_{L^2}
			\|\nabla\omega_{\neq}\|_{L^2}\\\le~&
			\frac1{64A}\|\nabla\omega_{\neq}\|_{L^2}^2
			+\frac CA
			\left(
			\|(\omega_{\mathrm{in}})_0\|_{L^\infty}
			+CK_{\neq}\right)^2
			\|\omega_{\neq}\|_{L^2}^2.
		\end{align*}
	Substituting this and
	\eqref{eq:joint-vorticity-density-bound}
	into \eqref{eq:joint-vorticity-L2-identity}, we obtain
	\begin{equation}\label{eq:joint-vorticity-basic-bound}
		\begin{aligned}
			\frac d{dt}\|\omega_{\neq}\|_{L^2}^2
			+\frac{63}{32A}
			\|\nabla\omega_{\neq}\|_{L^2}^2\le
			\frac CA\|n_{\neq}\|_{L^2}^2
			+\frac CA
			\left(
			\|(\omega_{\mathrm{in}})_0\|_{L^\infty}
			+CK_{\neq}
			\right)^2
			\|\omega_{\neq}\|_{L^2}^2.
		\end{aligned}
	\end{equation}
	
	\noindent{\bf (ii) The Fourier correction estimate.}
By \eqref{eq:joint-vorticity-equation}, the Fourier
differentiation identities and \eqref{eq:correction-shear},
differentiation of the vorticity correction in
\eqref{eq:full-energy} gives
		\begin{align}\label{eq:joint-vorticity-correction-identity}
		&\frac1{8\pi}\frac d{dt}
		\sum_{k\ne0}\int_{\mathbb R}
		\operatorname{sgn}(k)
		\arctan\bigl(A^{-\frac13}|k|^{-\frac13}\eta\bigr)
		|\widehat\omega(t,k,\eta)|^2d\eta\notag\\&+
		\frac1{8\pi}\sum_{k\ne0}\int_{\mathbb R}
		\frac{
			A^{-\frac13}|k|^{\frac23}
			|\widehat\omega(t,k,\eta)|^2
		}{
			1+A^{-\frac23}|k|^{-\frac23}\eta^2
		}d\eta\notag\\=~&
		-\frac2{8\pi A}
		\sum_{k\ne0}\int_{\mathbb R}
		\operatorname{sgn}(k)
		\arctan\bigl(A^{-\frac13}|k|^{-\frac13}\eta\bigr)
		(k^2+\eta^2)|\widehat\omega(t,k,\eta)|^2 d\eta\notag\\&-
		\frac2{8\pi A}\operatorname{Re}
		\sum_{k\ne0}\int_{\mathbb R}
		\operatorname{sgn}(k)
		\arctan\bigl(A^{-\frac13}|k|^{-\frac13}\eta\bigr)
		\mathrm ik\,
		\overline{\widehat\omega(t,k,\eta)}
		\widehat n(t,k,\eta)d\eta\\&+
		\frac2{8\pi A}\operatorname{Re}
		\sum_{k\ne0}\int_{\mathbb R}
		\operatorname{sgn}(k)
		\arctan\bigl(A^{-\frac13}|k|^{-\frac13}\eta\bigr)\notag\\&\times
		\overline{[\nabla\omega_{\neq}]^{\wedge}(t,k,\eta)}
		\cdot
		\left[
		P_{\neq}\bigl(v_{\neq}\omega_0
		+v_{\neq}\omega_{\neq}\bigr)
		\right]^{\wedge}(t,k,\eta)d\eta\notag\\&-
		\frac2{8\pi A}\operatorname{Re}
		\sum_{k\ne0}\int_{\mathbb R}
		\operatorname{sgn}(k)
		\arctan\bigl(A^{-\frac13}|k|^{-\frac13}\eta\bigr)
		\overline{\widehat\omega(t,k,\eta)}
		[v_{1,0}\partial_x\omega_{\neq}]^{\wedge}
		(t,k,\eta)d\eta\notag\\\notag&=:I_{1}+I_{2}+I_{3}+I_{4}.
	\end{align}
	We estimate the terms on the right-hand side in order.
	The bound $|\arctan z|\le\frac{\pi}{2}$ and
	Parseval's identity yield
	\begin{equation}\label{eq:joint-vorticity-diffusion-correction}
		\begin{aligned}
			|I_{1}|\leq\frac1{8A}\|\nabla\omega_{\neq}\|_{L^2}^2.
		\end{aligned}
	\end{equation}
	For the density term $I_{2}$, by Cauchy--Schwarz inequality,
	Parseval's identity and
	\eqref{eq:joint-multiplier-gradient}, we have
	\begin{equation}\label{eq:joint-vorticity-corrected-density-bound}
		\begin{aligned}
		|I_{2}|\le
			\frac1{8A}\|n_{\neq}\|_{L^2}
			\|\nabla\omega_{\neq}\|_{L^2}\le
			\frac1{64A}\|\nabla\omega_{\neq}\|_{L^2}^2
			+\frac CA\|n_{\neq}\|_{L^2}^2.
		\end{aligned}
	\end{equation}
	For the product $v_{\neq}\omega_0$, the same
	Cauchy--Schwarz estimate and
	\eqref{eq:joint-vorticity-mean-product} give
		\begin{align*}
			&\frac2{8\pi A}
			\left|
			\operatorname{Re}
			\sum_{k\ne0}\int_{\mathbb R}
			\operatorname{sgn}(k)
			\arctan\bigl(A^{-\frac13}|k|^{-\frac13}\eta\bigr)
			\overline{[\nabla\omega_{\neq}]^{\wedge}(t,k,\eta)}
			\cdot
			[P_{\neq}(v_{\neq}\omega_0)]^{\wedge}
			(t,k,\eta)d\eta
			\right|\\
			&\le
			\frac1{8A}
			\|P_{\neq}(v_{\neq}\omega_0)\|_{L^2}
			\|\nabla\omega_{\neq}\|_{L^2}\le
			\frac1{64A}\|\nabla\omega_{\neq}\|_{L^2}^2
			+\frac CA
			\left(
			\|(\omega_{\mathrm{in}})_0\|_{L^\infty}
			+CK_{\neq}
			\right)^2
			\|\omega_{\neq}\|_{L^2}^2.
		\end{align*}
	For the product $v_{\neq}\omega_{\neq}$,
	H\"older's inequality, \eqref{eq:Biot-bound} and interpolation inequality imply
		\begin{align*}
			\|P_{\neq}(v_{\neq}\omega_{\neq})\|_{L^2}\le
			\|v_{\neq}\|_{L^4}\|\omega_{\neq}\|_{L^4}\le
			C\|\omega_{\neq}\|_{L^2}^{\frac32}
			\left(
			\|\omega_{\neq}\|_{L^2}
			+\|\nabla\omega_{\neq}\|_{L^2}
			\right)^{\frac12}.
		\end{align*}
	Applying this and \eqref{eq:joint-multiplier-gradient},
	we obtain
	\begin{align*}
		&\frac2{8\pi A}
		\left|
		\operatorname{Re}
		\sum_{k\ne0}\int_{\mathbb R}
		\operatorname{sgn}(k)
		\arctan\bigl(A^{-\frac13}|k|^{-\frac13}\eta\bigr)
		\overline{[\nabla\omega_{\neq}]^{\wedge}(t,k,\eta)}
		\cdot
		[P_{\neq}(v_{\neq}\omega_{\neq})]^{\wedge}
		(t,k,\eta)d\eta
		\right|\\
		&\le
		\frac1{8A}
		\|P_{\neq}(v_{\neq}\omega_{\neq})\|_{L^2}
		\|\nabla\omega_{\neq}\|_{L^2}\le
		\frac CA\|\omega_{\neq}\|_{L^2}^2
		\|\nabla\omega_{\neq}\|_{L^2}
		+\frac CA\|\omega_{\neq}\|_{L^2}^{\frac32}
		\|\nabla\omega_{\neq}\|_{L^2}^{\frac32}\\
		&\le
		\frac1{32A}\|\nabla\omega_{\neq}\|_{L^2}^2
		+\frac CA
		\left(
		\|\omega_{\neq}\|_{L^2}^4
		+\|\omega_{\neq}\|_{L^2}^6
		\right).
	\end{align*}
	By \eqref{eq:joint-proof-bootstrap-consequences}, we arrive at
	\[
	\|\omega_{\neq}\|_{L^2}^4
	+\|\omega_{\neq}\|_{L^2}^6
	\le
	\left[
	\frac{8K_{\neq}}3+
	\left(\frac{8K_{\neq}}3\right)^2
	\right]\|\omega_{\neq}\|_{L^2}^2.
	\]
	Therefore,
		\begin{align*}
			&\frac2{8\pi A}
			\left|
			\operatorname{Re}
			\sum_{k\ne0}\int_{\mathbb R}
			\operatorname{sgn}(k)
			\arctan\bigl(A^{-\frac13}|k|^{-\frac13}\eta\bigr)
			\overline{[\nabla\omega_{\neq}]^{\wedge}(t,k,\eta)}
			\cdot
			[P_{\neq}(v_{\neq}\omega_{\neq})]^{\wedge}
			(t,k,\eta)d\eta
			\right|\\&\le
			\frac1{32A}\|\nabla\omega_{\neq}\|_{L^2}^2
			+\frac{C(K_{\neq})}A
			\|\omega_{\neq}\|_{L^2}^2.
		\end{align*}
	Collecting the above estimates, $|I_{3}|$ is bounded by
	\begin{equation}\label{I3 estimate}
	 \frac{3}{64A}\|\nabla\omega_{\neq}\|_{L^{2}}^{2}+\frac{C(K_{\neq},\|(\omega_{\mathrm{in}})_0\|_{L^{\infty}} )}{A}\|\omega_{\neq}\|_{L^{2}}^{2}
	\end{equation}
Finally, applying
	\eqref{eq:mean-transport-correction}
	with $f=\omega_{\neq}$, and using
	$\partial_yv_{1,0}=\omega_0$ and
	\eqref{eq:joint-proof-mean-bound}, we find
	\begin{equation}\label{eq:joint-vorticity-mean-transport-bound}
		\begin{aligned}
		|I_{4}|\le
			\frac C{8\pi A}
			\left(
			\|(\omega_{\mathrm{in}})_0\|_{L^\infty}
			+CK_{\neq}
			\right)
			\sum_{k\ne0}\int_{\mathbb R}
			A^{-\frac13}|k|^{\frac23}
			|\widehat\omega(t,k,\eta)|^2 d\eta.
		\end{aligned}
	\end{equation}
	Substituting
	\eqref{eq:joint-vorticity-diffusion-correction},
	\eqref{eq:joint-vorticity-corrected-density-bound},
	\eqref{I3 estimate}
	and \eqref{eq:joint-vorticity-mean-transport-bound}
	into \eqref{eq:joint-vorticity-correction-identity}
	gives
	\begin{equation}\label{eq:joint-vorticity-correction-bound}
		\begin{aligned}
			&\frac1{8\pi}\frac d{dt}
			\sum_{k\ne0}\int_{\mathbb R}
			\operatorname{sgn}(k)
			\arctan\bigl(A^{-\frac13}|k|^{-\frac13}\eta\bigr)
			|\widehat\omega(t,k,\eta)|^2 d\eta\\&\quad+
			\frac1{8\pi}\sum_{k\ne0}\int_{\mathbb R}
			\frac{
				A^{-\frac13}|k|^{\frac23}
				|\widehat\omega(t,k,\eta)|^2
			}{
				1+A^{-\frac23}|k|^{-\frac23}\eta^2
			}d\eta\\
			&\le
			\frac{3}{16A}\|\nabla\omega_{\neq}\|_{L^2}^2
			+\frac CA\|n_{\neq}\|_{L^2}^2
			+\frac{
				C(K_{\neq},
				\|(\omega_{\mathrm{in}})_0\|_{L^\infty(\mathbb R)})
			}A\|\omega_{\neq}\|_{L^2}^2\\
			&\quad+
			\frac C{8\pi A}
			\left(
			\|(\omega_{\mathrm{in}})_0\|_{L^\infty(\mathbb R)}
			+CK_{\neq}
			\right)
			\sum_{k\ne0}\int_{\mathbb R}
			A^{-\frac13}|k|^{\frac23}
			|\widehat\omega(t,k,\eta)|^2 d\eta.
		\end{aligned}
	\end{equation}
	Adding \eqref{eq:joint-vorticity-basic-bound}
	and \eqref{eq:joint-vorticity-correction-bound},
	and moving the gradient term on the right-hand
	side to the left, we obtain
	\begin{equation}\label{eq:joint-vorticity-energy-bound}
			\begin{aligned}
			&\frac d{dt}\left[
			\|\omega_{\neq}\|_{L^2}^2
			+\frac1{8\pi}
			\sum_{k\ne0}\int_{\mathbb R}
			\operatorname{sgn}(k)
			\arctan\bigl(A^{-\frac13}|k|^{-\frac13}\eta\bigr)
			|\widehat\omega(t,k,\eta)|^2 d\eta
			\right]\\
			&\quad+
			\frac{57}{32A}\|\nabla\omega_{\neq}\|_{L^2}^2
			+\frac1{8\pi}
			\sum_{k\ne0}\int_{\mathbb R}
			\frac{
				A^{-\frac13}|k|^{\frac23}
				|\widehat\omega(t,k,\eta)|^2
			}{
				1+A^{-\frac23}|k|^{-\frac23}\eta^2
			}d\eta\\
			&\le
			\frac CA\|n_{\neq}\|_{L^2}^2
			+\frac{
				C(K_{\neq},
				\|(\omega_{\mathrm{in}})_0\|_{L^\infty(\mathbb R)})
			}A\|\omega_{\neq}\|_{L^2}^2\\
			&\quad+
			\frac C{8\pi A}
			\left(
			\|(\omega_{\mathrm{in}})_0\|_{L^\infty(\mathbb R)}
			+CK_{\neq}
			\right)
			\sum_{k\ne0}\int_{\mathbb R}
			A^{-\frac13}|k|^{\frac23}
			|\widehat\omega(t,k,\eta)|^2 d\eta.
		\end{aligned}
	\end{equation}

	\noindent{\bf Step III. Combining the estimates.}
	Add \eqref{eq:joint-density-energy-bound}
	and \eqref{eq:joint-vorticity-energy-bound}.
	By \eqref{eq:full-energy}, the sum of their
	time derivatives is $E'(t)$.
	We next combine the Couette dissipation with
	a portion of the gradient dissipation.
	For the density, apply
	\eqref{eq:fourier-coercivity} with
	\[
	a=
	\frac1{
		8\pi(1+2K_\infty)^{\alpha-1}
		[1+2(1-\alpha)K_\infty]},
	\]
	multiply that inequality by
	$|\widehat n(t,k,\eta)|^2$,
	integrate over $\eta\in\mathbb R$ and sum over $k\ne0$.
	For the vorticity, use $a=\frac{1}{8\pi}$ and
	$|\widehat\omega(t,k,\eta)|^2$.
	Both coefficients belong to $(0,1]$.
	For the density gradient, we also use
	$\|\nabla n_{\neq}\|_{L^2}\le\|\nabla n\|_{L^2}$. Therefore,
	\begin{align}\label{eq:joint-after-frequency-coercivity}
		&E'(t)+\frac5{4A}
		\left(
		\|\nabla n\|_{L^2}^2
		+\|\nabla\omega_{\neq}\|_{L^2}^2
		\right)+\frac1{32\pi}
		\sum_{k\ne0}\int_{\mathbb R}
		A^{-\frac13}|k|^{\frac23}
		|\widehat\omega(t,k,\eta)|^2 d\eta
		\notag\\
		&+\frac1{
			32\pi(1+2K_\infty)^{\alpha-1}
			[1+2(1-\alpha)K_\infty]}
		\sum_{k\ne0}\int_{\mathbb R}
		A^{-\frac13}|k|^{\frac23}
		|\widehat n(t,k,\eta)|^2 d\eta
		\notag\\&\le
		\frac{C(\alpha,M)}A
		+\frac{C(\alpha,M,K_\infty)}A
		\|n_{\neq}\|_{L^2}^2+
		\frac{
			C(\alpha,K_{\neq},K_\infty,
			\|(\omega_{\mathrm{in}})_0\|_{L^\infty(\mathbb R)})
		}A\|\omega_{\neq}\|_{L^2}^2\\&+
		\frac CA
		\left(
		\|(\omega_{\mathrm{in}})_0\|_{L^\infty(\mathbb R)}
		+CK_{\neq}
		\right)
		\sum_{k\ne0}\int_{\mathbb R}
		A^{-\frac13}|k|^{\frac23}
		\notag\\\notag&\times
		\left(
		\frac{|\widehat n(t,k,\eta)|^2}{
			8\pi(1+2K_\infty)^{\alpha-1}
			[1+2(1-\alpha)K_\infty]}
		+\frac{|\widehat\omega(t,k,\eta)|^2}{8\pi}
		\right)d\eta.
	\end{align}
	It remains to absorb the three nonconstant terms
	on the right-hand side.
	Fix their constants and choose $A$ sufficiently
	large that
	\begin{align}
		&\frac CA
		\left(
		\|(\omega_{\mathrm{in}})_0\|_{L^\infty(\mathbb R)}
		+CK_{\neq}
		\right)
		\le\frac1{16},
		\label{eq:joint-absorb-mean-transport}\\
		&C(\alpha,M,K_\infty)A^{-\frac23}
		(1+2K_\infty)^{\alpha-1}
		[1+2(1-\alpha)K_\infty]
		\le\frac1{128\pi},
		\label{eq:joint-absorb-density}\\
		&C(\alpha,K_{\neq},K_\infty,
		\|(\omega_{\mathrm{in}})_0\|_{L^\infty(\mathbb R)})
		A^{-\frac23}
		\le\frac1{128\pi}.
		\label{eq:joint-absorb-vorticity}
	\end{align}
	First, \eqref{eq:joint-absorb-mean-transport}
	bounds the coefficients of the density and
	vorticity Fourier integrals in the last term
	of \eqref{eq:joint-after-frequency-coercivity}
	by
	\[
	\frac1{
		128\pi(1+2K_\infty)^{\alpha-1}
		[1+2(1-\alpha)K_\infty]}\quad\text{and}\quad
	\frac1{128\pi},
	\]
	respectively.
	Second, Parseval's identity and $|k|^{\frac23}\ge1$
	for $k\ne0$ give
	\begin{equation}\label{eq:joint-nonzero-norm-comparison}
		\frac1A\|f_{\neq}\|_{L^2}^2
		\le
		A^{-\frac23}
		\sum_{k\ne0}\int_{\mathbb R}
		A^{-\frac13}|k|^{\frac23}
		|\widehat f(k,\eta)|^2 d\eta,
		\qquad f=n\quad\text{or}\quad \omega.
	\end{equation}
	Applying \eqref{eq:joint-nonzero-norm-comparison}
	with $f=n$ and using
	\eqref{eq:joint-absorb-density}, we obtain
	\[
	\begin{aligned}
		&\frac{C(\alpha,M,K_\infty)}A
		\|n_{\neq}\|_{L^2}^2
		\le
		C(\alpha,M,K_\infty)A^{-\frac23}
		\sum_{k\ne0}\int_{\mathbb R}
		A^{-\frac13}|k|^{\frac23}
		|\widehat n(t,k,\eta)|^2 d\eta\\
		&\le
		\frac1{
			128\pi(1+2K_\infty)^{\alpha-1}
			[1+2(1-\alpha)K_\infty]}
		\sum_{k\ne0}\int_{\mathbb R}
		A^{-\frac13}|k|^{\frac23}
		|\widehat n(t,k,\eta)|^2 d\eta.
	\end{aligned}
	\]
	The same calculation with $f=\omega$ and
	\eqref{eq:joint-absorb-vorticity} gives
	\[
	\begin{aligned}
		\frac{
			C(\alpha,K_{\neq},K_\infty,
			\|(\omega_{\mathrm{in}})_0\|_{L^\infty(\mathbb R)})
		}A\|\omega_{\neq}\|_{L^2}^2\le
		\frac1{128\pi}
		\sum_{k\ne0}\int_{\mathbb R}
		A^{-\frac13}|k|^{\frac23}
		|\widehat\omega(t,k,\eta)|^2 d\eta.
	\end{aligned}
	\]
	Substituting these three bounds into \eqref{eq:joint-after-frequency-coercivity}, we conclude that
	\begin{align*}
		&E'(t)+\frac1A
		\left(
		\|\nabla n\|_{L^2}^2
		+\|\nabla\omega_{\neq}\|_{L^2}^2
		\right)+\frac1{64\pi}
		\sum_{k\ne0}\int_{\mathbb R}
		A^{-\frac13}|k|^{\frac23}
		|\widehat\omega(t,k,\eta)|^2 d\eta\\
		&+\frac1{
			64\pi(1+2K_\infty)^{\alpha-1}
			[1+2(1-\alpha)K_\infty]}
		\sum_{k\ne0}\int_{\mathbb R}
		A^{-\frac13}|k|^{\frac23}
		|\widehat n(t,k,\eta)|^2 d\eta
		\le\frac{C(\alpha,M)}A.
	\end{align*}
	This proves \eqref{eq:joint-differential}.
\end{proof}
Combining the preceding non-zero-mode dissipation with
one-dimensional interpolation for $n_0$ gives a uniform
bound for the complete energy.
\begin{corollary}\label{cor:joint-uniform}
Under the assumptions of Proposition \ref{prop:joint-energy},
after increasing its threshold for $A$ if necessary, we have
	\begin{equation}\label{eq:joint-ODE}
		E'(t)+\frac1A E(t)\le\frac{C(\alpha,M)}A,
		\qquad
		E(t)\le\max\{E(0),C(\alpha,M)\}.
	\end{equation}
	The constant on the right is independent of the bootstrap
	constants.
\end{corollary}
\begin{proof}
Suppose first that $M>0$.
By \eqref{eq:ordinary-mean-GN} and Young's inequality, we get
\[
4\pi\|n_0\|_{L^2(\mathbb R)}^2
\le
2\pi\|\partial_yn_0\|_{L^2(\mathbb R)}^2+C(M).
\]
This implies
\[
\begin{aligned}
	\|\nabla n\|_{L^2}^2=
	2\pi\|\partial_yn_0\|_{L^2(\mathbb R)}^2
	+\|\nabla n_{\neq}\|_{L^2}^2\ge
	4\pi\|n_0\|_{L^2(\mathbb R)}^2-C(M).
\end{aligned}
\]
Choose $A$ sufficiently large that
\[
A^{\frac23}\ge
128\pi(1+2K_\infty)^{\alpha-1}
[1+2(1-\alpha)K_\infty].
\]
Since
$(1+2K_\infty)^{\alpha-1}
[1+2(1-\alpha)K_\infty]\ge1$
and $|k|^{\frac23}\ge1$ for $k\ne0$,
Parseval's identity yields
\[
\begin{aligned}
	\frac1{
		64\pi(1+2K_\infty)^{\alpha-1}
		[1+2(1-\alpha)K_\infty]}
	&\sum_{k\ne0}\int_{\mathbb R}
	A^{-\frac13}|k|^{\frac23}
	|\widehat n(t,k,\eta)|^2 d\eta\ge\frac2A\|n_{\neq}\|_{L^2}^2,\\
	\frac1{64\pi}
	\sum_{k\ne0}\int_{\mathbb R}
	A^{-\frac13}|k|^{\frac23}&
	|\widehat\omega(t,k,\eta)|^2 d\eta
	\ge\frac2A\|\omega_{\neq}\|_{L^2}^2.
\end{aligned}
\]
Substituting these bounds into
\eqref{eq:joint-differential}, and using
$\|n\|_{L^2}^2=
2\pi\|n_0\|_{L^2(\mathbb R)}^2
+\|n_{\neq}\|_{L^2}^2$,
we obtain
\[
E'(t)+\frac2A
\left(
\|n\|_{L^2}^2+\|\omega_{\neq}\|_{L^2}^2
\right)
\le\frac{C(\alpha,M)}A.
\]
Hence it follows from \eqref{eq:full-energy-equivalence} that
\[E'(t)+\frac1A E(t)\le\frac{C(\alpha,M)}A.\]
Multiplying it by ${\rm e}^{\frac{t}{A}}$ and integrating over
$[0,t]$ gives
\[
\begin{aligned}
	E(t)\le
	{\rm e}^{-\frac{t}{A}}E(0)
	+\bigl(1-{\rm e}^{-\frac{t}{A}}\bigr)C(\alpha,M)\le
	\max\{E(0),C(\alpha,M)\}.
\end{aligned}
\]
If $M=0$, nonnegativity and mass conservation imply
$n=0$. Under the same choice of $A$,
\eqref{eq:joint-differential} therefore gives
\[
E'(t)+\frac2A\|\omega_{\neq}\|_{L^2}^2
\le\frac{C(\alpha,0)}A.
\]
The energy equivalence yields the same differential
inequality and bound for $E(t)$.
This proves \eqref{eq:joint-ODE}.
\end{proof}

\section{Proof of global existence}
\label{sec:full-global-proof}
The joint energy estimate supplies a uniform $L^2$
bound. We use this bound to obtain a uniform
$L^\infty$ estimate for the full density, improve the
bootstrap bounds, and apply the continuation criterion
in Proposition \ref{prop:H3-continuation}.
\subsection{A uniform $L^\infty$ bound for the density}
\begin{lemma}\label{lem:full-density-Moser}
Let $[0,T]$ be a compact interval of a classical
solution of \eqref{eq:rescaled-system}.
Assume that \eqref{eq:full-bootstrap} holds and that
$A$ satisfies the threshold in
Corollary \ref{cor:joint-uniform}. Then
	\begin{equation}\label{eq:full-density-Moser-bound}
		\sup_{0\le t\le T}\|n(t)\|_{L^\infty}
		\le
		C\bigl(
		\alpha,M,\|n_{\mathrm{in}}\|_{L^\infty},
		\|(\omega_{\mathrm{in}})_{\neq}\|_{L^2}
		\bigr).
	\end{equation}
	The constant is independent of $A$ and $T$.
\end{lemma}
\begin{proof}
	If $M=0$, nonnegativity and mass conservation give $n=0$.
	We therefore assume $M>0$.

\noindent{\bf Step I. The $L^r$ energy estimate.}
For $r\ge2$, testing $\eqref{eq:rescaled-system}_1$
with $rn^{r-1}$ gives
	\begin{align*}
		\frac{d}{dt}
		\int_{\mathbb T\times\mathbb R}n(t,x,y)^r dxdy~+&~
		\frac{r(r-1)}A
		\int_{\mathbb T\times\mathbb R}
		(1+n)^{\alpha-1}\bigl[1+(1-\alpha)n\bigr]
		n^{r-2}|\nabla n|^2 dxdy\\=~&
		\frac{r(r-1)}A
		\int_{\mathbb T\times\mathbb R}
		n^{r-1}(1+n)^\alpha\nabla c\cdot\nabla n dxdy.
	\end{align*}
	Young's inequality and
	Lemma \ref{lem:revised-coefficients} imply
	\begin{equation}\label{nr energy}
			\begin{aligned}
			\frac{d}{dt}
			\int_{\mathbb T\times\mathbb R}n^r dxdy&+
			\frac{r(r-1)}{2A}
			\int_{\mathbb T\times\mathbb R}
			(1+n)^{\alpha-1}\bigl[1+(1-\alpha)n\bigr]
			n^{r-2}|\nabla n|^2 dxdy\\&\le
			\frac{C_\alpha r^2}A
			\int_{\mathbb T\times\mathbb R}
			\bigl(n^r+n^{r+\alpha}\bigr)|\nabla c|^2 dxdy.
		\end{aligned}	
	\end{equation}
The two lower bounds in Lemma \ref{lem:revised-coefficients} yield
\[
\begin{aligned}
	2(1+n)^{\alpha-1}[1+(1-\alpha)n]
	\ge 1+(1-\alpha)(1+n)^\alpha\ge 1+(1-\alpha)n^\alpha.
\end{aligned}
\]
Hence
\[
(1+n)^{\alpha-1}[1+(1-\alpha)n]
\ge\frac12\bigl[1+(1-\alpha)n^\alpha\bigr].
\]
Using
\[
\begin{aligned}
	n^{r-2}|\nabla n|^2=\frac4{r^2}|\nabla n^{\frac{r}{2}}|^2,\qquad
	n^{r+\alpha-2}|\nabla n|^2
	=\frac4{(r+\alpha)^2}
	|\nabla n^{\frac{r+\alpha}{2}}|^2,
\end{aligned}
\]
we obtain, for $r\ge2$,
\begin{align}\label{result 1}
	&r(r-1)\int_{\mathbb T\times\mathbb R}
	(1+n)^{\alpha-1}[1+(1-\alpha)n]
	n^{r-2}|\nabla n|^2dxdy\nonumber\\
	&\ge \frac{2(r-1)}r
	\|\nabla n^{\frac{r}{2}}\|_{L^2}^2
	+\frac{2(1-\alpha)r(r-1)}{(r+\alpha)^2}
	\|\nabla n^{\frac{r+\alpha}{2}}\|_{L^2}^2\\\nonumber&\ge
	\|\nabla n^{\frac{r}{2}}\|_{L^2}^2
	+\frac{4(1-\alpha)}9
	\|\nabla n^{\frac{r+\alpha}{2}}\|_{L^2}^2.
\end{align}
The last inequality follows from
$\frac{r-1}r\ge\frac12$ and $r+\alpha\le\frac32r$ with $r\ge2$.
For $q=\frac{r}{2}$ or $q=\frac{r+\alpha}{2}$, using
H\"older's inequality,
\eqref{eq:full-elliptic-gradient} and
\eqref{eq:cylinder-interpolation}, we get
\[
\begin{aligned}
	&\frac{C_\alpha r^2}{A}
	\int_{\mathbb T\times\mathbb R}
	n^{2q}|\nabla c|^2 dxdy\\\le~&
	\frac{C_\alpha r^2}{A}
	\|n\|_{L^\infty(0,T;L^2)}^2
	\|n^q\|_{L^2}\|\nabla n^q\|_{L^2}
	+\frac{C_\alpha r^2}{A}
	\|n\|_{L^\infty(0,T;L^2)}^2
	\|n^q\|_{L^2}^2\\\le~&
	\frac{\varepsilon}{A}
	\|\nabla n^q\|_{L^2}^2
	+\frac{C(\alpha,\varepsilon)r^4}{A}
	\left(
	\|n\|_{L^\infty(0,T;L^2)}^4
	+\|n\|_{L^\infty(0,T;L^2)}^2
	\right)\|n^q\|_{L^2}^2,
\end{aligned}
\]
where we choose $\varepsilon=\frac14$ for $q=\frac{r}{2}$ and
$\varepsilon=\frac{1-\alpha}{9}$ for
$q=\frac{r+\alpha}{2}$.
Substituting \eqref{result 1} and the preceding
estimates into \eqref{nr energy}, and discarding the nonnegative term
$(4A)^{-1}\|\nabla n^{\frac{r}{2}}\|_{L^2}^2$, we obtain
	\begin{align}\label{eq:full-density-Moser-energy}
		&\frac d{dt}
		\int_{\mathbb T\times\mathbb R}n^r dxdy
		+\frac{1-\alpha}{9A}
		\|\nabla n^{\frac{r+\alpha}{2}}\|_{L^2}^2\nonumber\\\le~&
		\frac{C_\alpha r^4}{A}
		\left(
		\|n\|_{L^\infty(0,T;L^2)}^4
		+\|n\|_{L^\infty(0,T;L^2)}^2
		\right)
		\left(
		\int_{\mathbb T\times\mathbb R}n^r dxdy
		+\int_{\mathbb T\times\mathbb R}n^{r+\alpha}dxdy
		\right)\\\nonumber\le~&
		\frac{Cr^4}{A}
		\left(
		\int_{\mathbb T\times\mathbb R}n^r dxdy
		+\int_{\mathbb T\times\mathbb R}
		n^{r+\alpha}dxdy
		\right),
	\end{align}
	where we used \eqref{eq:full-energy-equivalence} and \eqref{eq:joint-ODE}. Moreover, the constant
 $C$ depends only on $\alpha$, $M$,
$\|n_{\mathrm{in}}\|_{L^2}$ and
$\|(\omega_{\mathrm{in}})_{\neq}\|_{L^2}$,
and is independent of $r$, $A$ and $T$.
	
For nonzero $f\in L^1\cap H^1(\mathbb T\times\mathbb R)$,
the cylinder Nash inequality gives
\begin{equation*}
	\|\nabla f\|_{L^2}^2
	\ge
	\frac{\|f\|_{L^2}^4}{C\|f\|_{L^1}^2}
	-C\|f\|_{L^1}^2.
\end{equation*}
Applying the above inequality
to $f=n^{\frac{r+\alpha}{2}}$, for each fixed $r\ge2$, we obtain
\[
\begin{aligned}
	&\|\nabla n^{\frac{r+\alpha}{2}}\|_{L^2}^2\ge
	\frac1C
	\frac{
		\left(
		\displaystyle\int_{\mathbb T\times\mathbb R}
		n^{r+\alpha}dxdy
		\right)^2
	}{
		\left(
		\displaystyle\int_{\mathbb T\times\mathbb R}
		n^{\frac{r+\alpha}{2}}dxdy
		\right)^2}-
	C\left(
	\int_{\mathbb T\times\mathbb R}
	n^{\frac{r+\alpha}{2}}dxdy
	\right)^2\\
	&\ge
	\frac1C
	\frac{
		\left(
		\displaystyle\int_{\mathbb T\times\mathbb R}
		n^{r+\alpha}dxdy
		\right)^2
	}{
		\left(
		\displaystyle\sup_{0\le s\le T}
		\int_{\mathbb T\times\mathbb R}
		n^{\frac{r+\alpha}{2}}dxdy
		\right)^2}-C\left(
	\sup_{0\le s\le T}
	\int_{\mathbb T\times\mathbb R}
	n^{\frac{r+\alpha}{2}}dxdy
	\right)^2.
\end{aligned}
\]
Substituting this bound into
\eqref{eq:full-density-Moser-energy} gives
\begin{equation}\label{eq:review-Moser-ODE}
	\begin{aligned}
		&\frac d{dt}
		\int_{\mathbb T\times\mathbb R}
		n^rdxdy\le
		\frac{Cr^4}{A}
		\left(
		\int_{\mathbb T\times\mathbb R}
		n^rdxdy
		+\int_{\mathbb T\times\mathbb R}
		n^{r+\alpha}dxdy
		\right)\\&+
		\frac{C_\alpha}{A}
		\left(
		\sup_{0\le s\le T}
		\int_{\mathbb T\times\mathbb R}
		n^{\frac{r+\alpha}{2}}dxdy
		\right)^2-
		\frac1{C_\alpha A}
		\frac{\left(
			\displaystyle\int_{\mathbb T\times\mathbb R}
			n^{r+\alpha}dxdy
			\right)^2}{
			\left(
			\displaystyle\sup_{0\le s\le T}
			\int_{\mathbb T\times\mathbb R}
			n^{\frac{r+\alpha}{2}}dxdy
			\right)^2}.
	\end{aligned}
\end{equation}
Noting that $n^r\le n^2$ on $\{n\le1\}$ and
$n^r\le n^{r+\alpha}$ on $\{n>1\}$, we have
\begin{equation}\label{eq:moser-moment-comparison}
	\int_{\mathbb T\times\mathbb R}n^rdxdy
	\le\|n\|_{L^\infty(0,T;L^2)}^2
	+\int_{\mathbb T\times\mathbb R}
	n^{r+\alpha}dxdy.
\end{equation}
Choose $C_1\ge2$ sufficiently large.
If
\[
\begin{aligned}
	\int_{\mathbb T\times\mathbb R}n^rdxdy
	>C_1r^4\biggl[1+\|n\|_{L^\infty(0,T;L^2)}^2+
	\left(
	\sup_{0\le s\le T}
	\int_{\mathbb T\times\mathbb R}
	n^{\frac{r+\alpha}{2}}dxdy
	\right)^2
	\biggr],
\end{aligned}
\]
then \eqref{eq:moser-moment-comparison} implies
\begin{equation}\label{eq:moser-above-threshold}
	\begin{aligned}
		\int_{\mathbb T\times\mathbb R}
		n^{r+\alpha}dxdy
		>\frac{C_1r^4}{2}\biggl[1+\|n\|_{L^\infty(0,T;L^2)}^2+
		\left(
		\sup_{0\le s\le T}
		\int_{\mathbb T\times\mathbb R}
		n^{\frac{r+\alpha}{2}}dxdy
		\right)^2
		\biggr].
	\end{aligned}
\end{equation}
In particular, $\int_{\mathbb T\times\mathbb R}n^rdxdy\le
2\int_{\mathbb T\times\mathbb R}n^{r+\alpha}dxdy$ and
\[
\begin{aligned}
	\left(
	\sup_{0\le s\le T}
	\int_{\mathbb T\times\mathbb R}
	n^{\frac{r+\alpha}{2}}dxdy
	\right)^2\le
	\frac2{C_1r^4}
	\int_{\mathbb T\times\mathbb R}
	n^{r+\alpha}dxdy.
\end{aligned}
\]
Moreover, \eqref{eq:moser-above-threshold} gives
\[
\frac{
	\left(
	\displaystyle\int_{\mathbb T\times\mathbb R}
	n^{r+\alpha}dxdy
	\right)^2
}{\left(
	\displaystyle\sup_{0\le s\le T}
	\int_{\mathbb T\times\mathbb R}
	n(s,x,y)^{\frac{r+\alpha}{2}}dxdy
	\right)^2}\ge
\frac{C_1r^4}{2}
\int_{\mathbb T\times\mathbb R}
n^{r+\alpha}dxdy.
\]
Consequently, \eqref{eq:review-Moser-ODE} yields
\[
\begin{aligned}
	\frac d{dt}
	\int_{\mathbb T\times\mathbb R}n^rdxdy
	&\le
	\frac1A
	\left(
	3Cr^4+\frac{2C_\alpha}{C_1r^4}
	-\frac{C_1r^4}{2C_\alpha}
	\right)
	\int_{\mathbb T\times\mathbb R}
	n^{r+\alpha}dxdy\\
	&\le
	\frac{r^4}{A}
	\left(
	3C+C_\alpha-\frac{C_1}{2C_\alpha}
	\right)
	\int_{\mathbb T\times\mathbb R}
	n^{r+\alpha}dxdy.
\end{aligned}
\]
Choose $C_1>2C_\alpha(3C+C_\alpha)$.
The preceding estimate then gives
\[
\frac d{dt}
\int_{\mathbb T\times\mathbb R}n^rdxdy<0.
\]
Since the integral is absolutely continuous in time,
it is decreasing on every interval where it exceeds
that threshold. Consequently, its value is bounded
by the maximum of its initial value and the threshold.
Writing $C$ for the chosen constant $C_1$, we obtain
\begin{equation}\label{eq:full-density-Moser-recursion}
	\begin{aligned}
		&\sup_{0\le t\le T}
		\int_{\mathbb T\times\mathbb R}
		n^rdxdy\\\le~&
		\max\biggl\{
		\int_{\mathbb T\times\mathbb R}
		n_{\mathrm{in}}^rdxdy, 
		Cr^4\biggl[
		1+\|n\|_{L^\infty(0,T;L^2)}^2+
		\left(
		\sup_{0\le s\le T}
		\int_{\mathbb T\times\mathbb R}
		n^{\frac{r+\alpha}{2}}dxdy
		\right)^2
		\biggr]\biggr\}.
	\end{aligned}
\end{equation}

\noindent{\bf Step II. The Moser iteration.}
Set
$r_j=\alpha+(2-\alpha)2^j$ for $j\in\{0,1,\ldots\}$.
Then
$r_0=2$ and $\frac{r_j+\alpha}{2}=r_{j-1}$ with $j\ge1$.
For every $j\ge0$, mass conservation at the initial
time gives
\[
\begin{aligned}
	\int_{\mathbb T\times\mathbb R}
	n_{\mathrm{in}}(x,y)^{r_j}dxdy\le
	M\|n_{\mathrm{in}}\|_{L^\infty}^{r_j-1}\le
	\max\bigl\{
	1,M,\|n_{\mathrm{in}}\|_{L^\infty}
	\bigr\}^{r_j}.
\end{aligned}
\]
Applying \eqref{eq:full-density-Moser-recursion}
with $r=r_j$ and taking $r_j$-th roots, for $j\geq 1$, we obtain
\begin{equation}
	\begin{aligned}
		&\max\biggl\{
		1,M,\|n_{\mathrm{in}}\|_{L^\infty},\|n\|_{L^{\infty}(0,T;L^{2})}, \sup_{0\le t\le T}\|n(t)\|_{L^{r_j}}
		\biggr\}\\\le~&
		(Cr_j^4)^{\frac{1}{r_j}}
		\max\biggl\{
		1,M,\|n_{\mathrm{in}}\|_{L^\infty},
		\|n\|_{L^\infty(0,T;L^2)},
		\sup_{0\le t\le T}\|n(t)\|_{L^{r_{j-1}}}
		\biggr\}^{\frac{2r_{j-1}}{r_j}},
		\label{eq:moser-norm-recursion}
	\end{aligned}
\end{equation}
The products of the iteration exponents satisfy
$\prod_{\ell=1}^j\frac{2r_{\ell-1}}{r_\ell}
=\frac{2^{j+1}}{r_j}
\le\frac2{2-\alpha}$.
Similarly, for $1\le i\le j$,
\[
\prod_{\ell=i+1}^j
\frac{2r_{\ell-1}}{r_\ell}
=
\frac{2^{j-i}r_i}{r_j}
\le\frac2{2-\alpha},
\]
where the product equals one when $i=j$.
Taking logarithms in
\eqref{eq:moser-norm-recursion} and iterating,
we therefore obtain
\[
\begin{aligned}
	&\log\max\biggl\{
	1,M,\|n_{\mathrm{in}}\|_{L^\infty},
	\|n\|_{L^\infty(0,T;L^2)},\;
	\sup_{0\le t\le T}\|n(t)\|_{L^{r_j}}
	\biggr\}\\\le~&
	\frac2{2-\alpha}
	\log\max\bigl\{
	1,M,\|n_{\mathrm{in}}\|_{L^\infty},
	\|n\|_{L^\infty(0,T;L^2)}
	\bigr\}+
	\frac2{2-\alpha}
	\sum_{i=1}^j
	\frac{\log C+4\log r_i}{r_i}.
\end{aligned}
\]
Due to $(2-\alpha)2^i\le r_i\le2^{i+1}$,
we have
\[
\sum_{i=1}^\infty
\frac{\log C+4\log r_i}{r_i}
\le
\frac1{2-\alpha}
\sum_{i=1}^\infty
\frac{\log C+4(i+1)\log2}{2^i}
<\infty.
\]
The preceding logarithmic estimate and the bound
for the series imply
\[
\begin{aligned}
	\sup_{0\le t\le T}\|n(t)\|_{L^{r_j}}
	\le
	\exp\left(
	\frac{2(\log C+12\log2)}{(2-\alpha)^2}
	\right)
	\max\bigl\{
	1,M,\|n_{\mathrm{in}}\|_{L^\infty},
	\|n\|_{L^\infty(0,T;L^2)}
	\bigr\}^{\frac2{2-\alpha}}.
\end{aligned}
\]
Since the right-hand side is independent of $j$, there holds
\[
\sup_{j\ge0}\sup_{0\le t\le T}
\|n(t)\|_{L^{r_j}}
\le
C\left(
\alpha,M,\|n_{\mathrm{in}}\|_{L^\infty},
\|n\|_{L^\infty(0,T;L^2)}
\right).
\]
Letting $j\to\infty$ yields
\[
\|n\|_{L^\infty(0,T;L^\infty)}
\le
C\left(
\alpha,M,\|n_{\mathrm{in}}\|_{L^\infty},
\|n\|_{L^\infty(0,T;L^2)}
\right).
\]
By \eqref{eq:full-energy-equivalence}, \eqref{eq:joint-ODE} and
$\|n_{\mathrm{in}}\|_{L^2}^2
\le M\|n_{\mathrm{in}}\|_{L^\infty}$, we have
\[
\begin{aligned}
	\|n\|_{L^\infty(0,T;L^2)}^2&\le
	\frac43\max\{E(0),C(\alpha,M)\}\\&\le
	\frac43\max\left\{
	\frac54\left(
	M\|n_{\mathrm{in}}\|_{L^\infty}
	+\|(\omega_{\mathrm{in}})_{\neq}\|_{L^2}^2
	\right),C(\alpha,M)
	\right\}.
\end{aligned}
\]
Substituting this bound into the preceding
$L^\infty$ estimate proves
\eqref{eq:full-density-Moser-bound}.
\end{proof}

\subsection{Choice of constants and bootstrap improvements}
Let $T_{\max}$ be the maximal existence time of
the solution in Proposition \ref{prop:H3-continuation}.
We choose the constants introduced in
Section~\ref{sec:main-ideas} as follows:
\begin{equation}\label{eq:choose-levels}
	\begin{aligned}
		K_{\neq}
		&=4\left[
		1+\frac54\left(
		\|n_{\mathrm{in}}\|_{L^2}^2
		+\|(\omega_{\mathrm{in}})_{\neq}\|_{L^2}^2
		\right)+C(\alpha,M)
		\right],\\
		K_\infty
		&=2\left[
		1+\|n_{\mathrm{in}}\|_{L^\infty}
		+C\left(
		\alpha,M,\|n_{\mathrm{in}}\|_{L^\infty},
		\|(\omega_{\mathrm{in}})_{\neq}\|_{L^2}
		\right)
		\right].
	\end{aligned}
\end{equation}
Here the constants are those in
\eqref{eq:joint-ODE} and
\eqref{eq:full-density-Moser-bound}, respectively.
Both bootstrap constants are independent of $A$ and $T$.
Choose $A_0\ge1$ so that the absorption conditions
in Proposition \ref{prop:joint-energy} and
Corollary \ref{cor:joint-uniform} hold for these
fixed constants.
Since
\[
\begin{aligned}
	\|n_{\mathrm{in}}\|_{L^\infty}\le C\|n_{\mathrm{in}}\|_{H^3},\quad
	\|(\omega_{\mathrm{in}})_{\neq}\|_{L^2}
	+\|(\omega_{\mathrm{in}})_{0}\|_{L^\infty(\mathbb R)}\le C\|v_{\mathrm{in}}\|_{H^3},
\end{aligned}
\]
$A_0$ depends only on the quantities stated in
Theorem \ref{thm:full-global}.

Let $A\ge A_0$ and assume that
\eqref{eq:full-bootstrap} holds on $[0,T]$.
By \eqref{eq:full-energy-equivalence}, \eqref{eq:joint-ODE} and
\eqref{eq:choose-levels}, we have
\[
E(t)\le\max\{E(0),C(\alpha,M)\}
\le\frac{K_{\neq}}4,
\qquad 0\le t\le T.
\]
Consequently,
\begin{equation}\label{eq:improved-energy}
	\begin{aligned}
		\sup_{0\le t\le T}E(t)\le\frac{K_{\neq}}4,\qquad
		\sup_{0\le t\le T}
		\left(
		\|n(t)\|_{L^2}^2
		+\|\omega_{\neq}(t)\|_{L^2}^2
		\right)\le\frac{K_{\neq}}3.
	\end{aligned}
\end{equation}
Lemma \ref{lem:full-density-Moser} and
\eqref{eq:choose-levels} also give
\begin{equation}\label{eq:improved-density}
	\begin{aligned}
		\sup_{0\le t\le T}\|n(t)\|_{L^\infty}\le
		C\left(
		\alpha,M,\|n_{\mathrm{in}}\|_{L^\infty},
		\|(\omega_{\mathrm{in}})_{\neq}\|_{L^2}
		\right)\le\frac{K_\infty}{2}.
	\end{aligned}
\end{equation}
The bootstrap inequalities hold strictly at $t=0$.
By continuity, they hold on an initial time interval.
A first equality in either bootstrap bound before
the maximal existence time would contradict
\eqref{eq:improved-energy} or
\eqref{eq:improved-density}.
Thus the bootstrap bounds hold on every compact
interval of the maximal solution.
\begin{proof}[Proof of Theorem \ref{thm:full-global}]
	For $A\ge A_0$, the bootstrap argument above gives
	\eqref{eq:improved-energy} and \eqref{eq:improved-density}
	on every compact interval of the maximal solution.
	Hence the density has bounded $L^2$ and $L^\infty$
	norms throughout its existence interval.
	Proposition~\ref{prop:H3-continuation} implies
	$T_{\max}=\infty$.
	
	The uniform bounds for $n$ and $\omega_{\neq}$ follow
	from \eqref{eq:improved-energy} and
	\eqref{eq:improved-density}.
	Lemma \ref{lem:mean-vorticity-uniform} gives the
	uniform bound for $\omega_0$.
	Nonnegativity and mass conservation follow from
	Lemma \ref{lem:revised-mass}.
	Finally, reversing the time rescaling proves the
	assertions for \eqref{eq:full-model-revised}.
\end{proof}

\appendix
\section{Local existence and continuation}
\label{sec:local-continuation}
We prove Proposition \ref{prop:H3-continuation} using
the local quasilinear parabolic theory of
Amann \cite{Amann2005,Amann2017} and the scalar parabolic
regularity estimates in \cite{LSU}.
The constants in this appendix may depend on $A$
and on the finite time interval under consideration.
\begin{proof}[Proof of Proposition \ref{prop:H3-continuation}]
\noindent\textbf{1. Local existence.}
The perturbation velocity satisfies
\[
\partial_tv+y\partial_xv+(v_2,0)
+\frac1A v\cdot\nabla v+\frac1A\nabla p
=\frac1A\Delta v+\frac1Ane_2,
\qquad \nabla\cdot v=0.
\]
We make the coordinate change $(x,y)\mapsto(x+ty,y)$.
For a scalar function $f$,
\[
\partial_t\bigl[f(t,x+ty,y)\bigr]
=(\partial_tf+y\partial_xf)(t,x+ty,y).
\]
We compose $n$, $c$ and $p$ with this map and replace
the velocity by $(v_1-tv_2,v_2)(t,x+ty,y)$.
The latter transformation preserves incompressibility, since
\[
\begin{aligned}
	\partial_x\bigl[(v_1-tv_2)(t,x+ty,y)\bigr]
	+\partial_y\bigl[v_2(t,x+ty,y)\bigr]=(\partial_xv_1+\partial_yv_2)(t,x+ty,y)=0.
\end{aligned}
\]
We denote the transformed functions by the same symbols.
The original spatial derivatives become $\partial_x$ and
$\partial_y-t\partial_x$, and therefore
\begin{equation*}
	\partial_x^2+(\partial_y-t\partial_x)^2
	=\nabla\cdot\left[
	\begin{pmatrix}1+t^2&-t\\-t&1\end{pmatrix}\nabla
	\right].
\end{equation*}
The matrix has determinant one and trace $2+t^2$. Thus, for
$0\le t\le T$,
\begin{equation}\label{eq:local-transformed-ellipticity}
	\frac{k^2+\eta^2}{2+T^2}
	\le k^2+(\eta-tk)^2
	\le(2+T^2)(k^2+\eta^2).
\end{equation}
In particular, the transformed diffusion is uniformly elliptic
on finite time intervals. The coordinate and velocity
transformations, and their inverses, are bounded on $H^j$
and $W^{j,6}$, with bounds uniform on such intervals.
The coordinate change also preserves Lebesgue measure.
The transformed equations are
\begin{equation}\label{eq:local-transformed-system}
	\left\{
	\begin{aligned}
		~&\partial_tn=\frac1A\nabla\cdot\left[
		\begin{pmatrix}1+t^2&-t\\-t&1\end{pmatrix}
		(1+n)^{\alpha-1}[1+(1-\alpha)n]\nabla n
		\right]\\
		&\quad\quad\quad-\frac1A\nabla\cdot\left[
		\begin{pmatrix}1+t^2&-t\\-t&1\end{pmatrix}
		n(1+n)^\alpha\nabla c+vn\right],\\~&
		\bigl[1-\partial_x^2-(\partial_y-t\partial_x)^2\bigr]c
		=n,\\~&
		\partial_tv-\frac1A
		\bigl[\partial_x^2+(\partial_y-t\partial_x)^2\bigr]v=-2(v_2,0)-\frac1A v\cdot\nabla v\\
		&\quad\quad\quad+\frac1A n(-t,1)-\frac1A
		\begin{pmatrix}1+t^2&-t\\-t&1\end{pmatrix}\nabla p,
		\\~&
		\nabla\cdot v=0.
	\end{aligned}
	\right.
\end{equation}
The term $-2(v_2,0)$ combines the term $-(v_2,0)$ in
the original equation with the term $-(v_2,0)$ obtained
by differentiating $v_1-tv_2$ in time.
At $t=0$, both transformations are the identity.
To eliminate the pressure, we apply the Fourier multiplier
with symbol
\begin{equation*}
	I-\frac{
		\begin{pmatrix}1+t^2&-t\\-t&1\end{pmatrix}
		\binom{k}{\eta}(k,\eta)
	}{k^2+(\eta-tk)^2},
	\qquad (k,\eta)\ne(0,0),
\end{equation*}
and value $I$ at the origin.
Incompressibility gives
$(k,\eta)\widehat v=0$ and $(k,\eta)\partial_t\widehat v=0$,
so this multiplier leaves $v$ and $\partial_tv$ unchanged.
Moreover,
\[
(k,\eta)
\begin{pmatrix}1+t^2&-t\\-t&1\end{pmatrix}
\binom{k}{\eta}=k^2+(\eta-tk)^2.
\]
Consequently, it annihilates the transformed pressure
gradient, whose Fourier transform is
$\mathrm i
\begin{pmatrix}1+t^2&-t\\-t&1\end{pmatrix}
\binom{k}{\eta}\widehat p$.
It commutes with the transformed Laplacian because the
latter has a scalar Fourier symbol. The multiplier is the
usual Leray projection conjugated by the coordinate and
velocity transformations. It is therefore bounded on $L^6$
and $H^j$, uniformly on finite time intervals.

We apply the local parabolic theory with time integrability
exponent six and spatial spaces $L^6$ and $W^{2,6}$.
The corresponding initial-data space is $W^{\frac53,6}$, and
\begin{equation}\label{eq:local-initial-embedding}
	H^3(\mathbb T\times\mathbb R)
	\hookrightarrow W^{\frac53,6}(\mathbb T\times\mathbb R)
	\hookrightarrow W^{1,\infty}(\mathbb T\times\mathbb R).
\end{equation}
For the velocity, we use the divergence-free subspaces.
Since $n_{\mathrm{in}}\ge0$, a sufficiently small bounded
neighborhood of $n_{\mathrm{in}}$ in $W^{\frac53,6}$ satisfies
$n>-1/2$. On this neighborhood the coefficient functions
are smooth with bounded derivatives, and
$(1+n)^{\alpha-1}[1+(1-\alpha)n]\ge\frac12$.
Together with \eqref{eq:local-transformed-ellipticity}, this
gives uniform ellipticity of the density equation.
Expanding its diffusion and chemotactic terms gives
\[
\begin{aligned}
	&\nabla\cdot\left[
	\begin{pmatrix}1+t^2&-t\\-t&1\end{pmatrix}
	(1+n)^{\alpha-1}[1+(1-\alpha)n]\nabla n\right]=(1+n)^{\alpha-1}[1+(1-\alpha)n]\\&\times
	\bigl[\partial_x^2+(\partial_y-t\partial_x)^2\bigr]n+\alpha(1-\alpha)n(1+n)^{\alpha-2}
	\left(|\partial_xn|^2
	+|(\partial_y-t\partial_x)n|^2\right),\\[1mm]
	&\nabla\cdot\left[
	\begin{pmatrix}1+t^2&-t\\-t&1\end{pmatrix}
	n(1+n)^\alpha\nabla c\right]\\
	&=\bigl[(1+n)^\alpha
	+\alpha n(1+n)^{\alpha-1}\bigr]
	\left[\partial_xn\,\partial_xc
	+(\partial_y-t\partial_x)n\,(\partial_y-t\partial_x)c\right]+n(1+n)^\alpha(c-n).
\end{aligned}
\]
Here the last identity uses the chemical equation in
\eqref{eq:local-transformed-system}. The elliptic estimates
on each finite interval give
\[
\|c\|_{W^{2,6}}\le C_T\|n\|_{L^6},
\qquad
\|\nabla c\|_{L^\infty}\le C_T\|n\|_{L^\infty}.
\]
The product estimates
\[
\begin{aligned}
	\||\nabla n|^2\|_{L^6}
	\le\|\nabla n\|_{L^\infty}\|\nabla n\|_{L^6},\qquad
	\|v\cdot\nabla n\|_{L^6}
	\le\|v\|_{L^\infty}\|\nabla n\|_{L^6},
\end{aligned}
\]
$$\|v\cdot\nabla v\|_{L^6}
\le\|v\|_{L^\infty}\|\nabla v\|_{L^6} $$
and the same estimates for differences show that the
terms other than the principal diffusion are locally
Lipschitz from $W^{\frac53,6}\times W^{\frac53,6}$ into
$L^6\times L^6$. The principal density operator also
depends locally Lipschitz continuously on $n$ as an
operator from $W^{2,6}$ into $L^6$.
Freezing its coefficient at $n_{\mathrm{in}}$ gives
\[
-\frac1A(1+n_{\mathrm{in}})^{\alpha-1}
[1+(1-\alpha)n_{\mathrm{in}}]
\bigl[\partial_x^2+(\partial_y-t\partial_x)^2\bigr].
\]
Its coefficients have bounded spatial H\"older derivatives
by $n_{\mathrm{in}}\in H^3$, are smooth in time, and are
uniformly elliptic. The linear maximal regularity theorem
\cite[Theorem 1.23]{Amann2017} applies on the flat cylinder.
It also applies to the transformed heat operator, which
preserves the divergence-free subspace.
For $n$ sufficiently close to $n_{\mathrm{in}}$ in
$C([0,T_*];W^{\frac53,6})$, the coefficient difference satisfies
\[
\begin{aligned}
	&\left\|\left\{
	(1+n)^{\alpha-1}[1+(1-\alpha)n]
	-(1+n_{\mathrm{in}})^{\alpha-1}
	[1+(1-\alpha)n_{\mathrm{in}}]
	\right\}
	\bigl[\partial_x^2+(\partial_y-t\partial_x)^2\bigr]h
	\right\|_{L^6}\\
	&\le C\|n-n_{\mathrm{in}}\|_{L^\infty}
	\|h\|_{W^{2,6}},\qquad h\in W^{2,6}.
\end{aligned}
\]
The principal density operator is therefore a small
perturbation of the frozen operator from $W^{2,6}$ to $L^6$,
uniformly in time, and retains maximal regularity.
The remaining terms belong to $L^\infty(0,T_*;L^6)$
for continuous $W^{\frac53,6}$ trajectories.
Thus the hypotheses of \cite[Theorem 2.1]{Amann2005}
hold with time exponent six and forcing exponent infinity.
This theorem yields a unique local
solution satisfying
\begin{equation}\label{eq:local-maximal-regularity}
	n,v\in W^{1,6}(0,T_*;L^6)
	\cap L^6(0,T_*;W^{2,6})
	\cap C([0,T_*];W^{\frac53,6}).
\end{equation}
We choose $T_*$ strictly inside this existence interval.
Then \eqref{eq:local-initial-embedding} bounds $n$, $v$
and their first spatial derivatives on $[0,T_*]$.

\noindent\textbf{2. The $H^3$ estimate and local regularity.}
We first prove an estimate for smooth solutions in the
original coordinates. Suppose that $n>-\frac12$ and that
$n$, $\nabla n$ and $\nabla v$ are bounded on the interval
under consideration. In this estimate, $C$ depends on
$\alpha$, the fixed value of $A$, and these bounds, but
not on higher Sobolev norms. For $|\ell|\le3$, the Sobolev and commutator estimates give
\begin{equation}\label{eq:local-H3-product-estimates}
	\begin{aligned}
		&\|(1+n)^{\alpha-1}[1+(1-\alpha)n]-1\|_{H^3}
		+\|n(1+n)^\alpha\|_{H^3}\le C\|n\|_{H^3},\\
		&\left\|\left[\partial^\ell,
		(1+n)^{\alpha-1}[1+(1-\alpha)n]\right]\nabla n
		\right\|_{L^2}\le C\|n\|_{H^3},\\
		&\|[\partial^\ell,v\cdot\nabla]n\|_{L^2}
		\le C\left(
		\|\nabla v\|_{L^\infty}\|n\|_{H^3}
		+\|\nabla n\|_{L^\infty}\|v\|_{H^3}\right),\\
		&\|[\partial^\ell,v\cdot\nabla]v\|_{L^2}
		\le C\|\nabla v\|_{L^\infty}\|v\|_{H^3}.
	\end{aligned}
\end{equation}
The subtraction of one in the diffusion-coefficient
estimate is necessary because a nonzero constant does
not belong to $L^2(\mathbb T\times\mathbb R)$.
The elliptic estimates and the Sobolev product estimate
further imply
\begin{equation}\label{eq:local-H3-chemotaxis}
	\begin{aligned}
		\|n(1+n)^\alpha\nabla c\|_{H^3}
		&\le C\|n(1+n)^\alpha\|_{L^\infty}\|\nabla c\|_{H^3}
		+C\|n(1+n)^\alpha\|_{H^3}\|\nabla c\|_{L^\infty}\\
		&\le C\|n\|_{H^3},
	\end{aligned}
\end{equation}
where
$\|\nabla c\|_{H^3}\le C\|n\|_{H^2}$ and
$\|\nabla c\|_{L^\infty}\le C\|n\|_{L^\infty}$.
Apply $\partial^\ell$ to the density and velocity
equations, test with $2\partial^\ell n$ and
$2\partial^\ell v$, and sum over $|\ell|\le3$.
Incompressibility eliminates the leading velocity
transport and pressure terms. For the Couette transport, we have
\[
[\partial^\ell,y\partial_x]f
=\ell_2\partial_x^{\ell_1+1}\partial_y^{\ell_2-1}f
\quad\text{if}\quad\ell_2\ge1,
\]
and the commutator vanishes if $\ell_2=0$.
These terms are therefore bounded by the $H^3$ energy.
Using \eqref{eq:local-H3-product-estimates},
\eqref{eq:local-H3-chemotaxis} and the diffusion lower
bound $\frac12$, we obtain
\[
\begin{aligned}
	&\frac d{dt}\bigl(\|n\|_{H^3}^2+\|v\|_{H^3}^2\bigr)
	+\frac1A\|\nabla n\|_{H^3}^2
	+\frac2A\|\nabla v\|_{H^3}^2\\\le~& C\bigl(\|n\|_{H^3}^2+\|v\|_{H^3}^2\bigr)
	+\frac CA\|n\|_{H^3}\|\nabla n\|_{H^3}.
\end{aligned}
\]
Young's inequality bounds the last term by
$(2A)^{-1}\|\nabla n\|_{H^3}^2+CA^{-1}\|n\|_{H^3}^2$.
Consequently,
\begin{equation}\label{eq:continuation-H3-energy}
	\begin{aligned}
		\frac d{dt}\bigl(\|n\|_{H^3}^2+\|v\|_{H^3}^2\bigr)
		+\frac1{2A}\bigl(
		\|\nabla n\|_{H^3}^2+\|\nabla v\|_{H^3}^2\bigr)\le C\bigl(\|n\|_{H^3}^2+\|v\|_{H^3}^2\bigr).
	\end{aligned}
\end{equation}
We now approximate the initial data in $H^3$ by smooth
data, preserving nonnegativity and incompressibility.
By \cite[Theorem 3.1]{Amann2005}, the corresponding
solutions exist on a common interval $[0,T_*]$ and
have uniform bounds in \eqref{eq:local-maximal-regularity}.
After decreasing $T_*$ if necessary, their densities
remain greater than $-\frac12$. The embedding
\eqref{eq:local-initial-embedding} gives uniform bounds
for their first spatial derivatives.
Applying \eqref{eq:continuation-H3-energy} and
Gronwall's inequality yields
\begin{equation}\label{eq:local-H3-uniform}
	\begin{aligned}
		\sup_{0\le t\le T_*}
		\bigl(\|n(t)\|_{H^3}^2+\|v(t)\|_{H^3}^2\bigr)+\frac1A\int_0^{T_*}
		\bigl(\|n(t)\|_{H^4}^2+\|v(t)\|_{H^4}^2\bigr)dt
		\le C,
	\end{aligned}
\end{equation}
with $C$ independent of the approximation. 
Continuous dependence identifies the limit with the
solution in \eqref{eq:local-maximal-regularity}.
Weak compactness and \eqref{eq:local-H3-uniform} give
$n,v\in L^\infty(0,T_*;H^3)\cap L^2(0,T_*;H^4)$.
The transformed equations, after eliminating the
pressure, then give
$\partial_tn,\partial_tv\in L^2(0,T_*;H^2)$.
The time-continuity theorem
\begin{equation}\label{eq:local-H3-time-continuity}
	L^2(0,T_*;H^4)\cap W^{1,2}(0,T_*;H^2)
	\hookrightarrow C([0,T_*];H^3)
\end{equation}
therefore proves
\[
n,v\in C([0,T_*];H^3)\cap L^2(0,T_*;H^4).
\]
Returning to the original variables preserves this
regularity. Parabolic regularity gives a classical
solution for positive times, and
Lemma \ref{lem:revised-mass} gives nonnegativity and
mass conservation.

\noindent\textbf{3. Finite-time estimates.}
Suppose that $T_{\max}<\infty$ and
\eqref{eq:continuation-density-assumption} holds, that is,
\[
\sup_{0\le t<T_{\max}}
\bigl(\|n(t)\|_{L^2}+\|n(t)\|_{L^\infty}\bigr)<\infty.
\]
In the remainder of the proof, constants may depend
on $\alpha$, $A$, $T_{\max}$, the initial data and these
density bounds, but are independent of $T<T_{\max}$.

\noindent{\bf (a) Velocity and vorticity energy estimates.}
Testing the full vorticity equation with $2\omega$ gives
\[
\frac d{dt}\|\omega\|_{L^2}^2
+\frac2A\|\nabla\omega\|_{L^2}^2
=\frac2A\int_{\mathbb T\times\mathbb R}
n_{\neq}\partial_x\omega dxdy.
\]
Here incompressibility eliminates the transport term.
Young's inequality yields
\begin{equation}\label{eq:continuation-vorticity}
	\frac d{dt}\|\omega\|_{L^2}^2
	+\frac1A\|\nabla\omega\|_{L^2}^2
	\le\frac1A\|n_{\neq}\|_{L^2}^2.
\end{equation}
Testing the perturbation velocity equation with $2v$
similarly gives
\[
\begin{aligned}
	\frac d{dt}\|v\|_{L^2}^2
	+\frac2A\|\nabla v\|_{L^2}^2&=-2\int_{\mathbb T\times\mathbb R}v_1v_2dxdy
	+\frac2A\int_{\mathbb T\times\mathbb R}nv_2 dxdy\\
	&\le\left(1+\frac1A\right)\|v\|_{L^2}^2
	+\frac1A\|n\|_{L^2}^2.
\end{aligned}
\]
For divergence-free $v$, Parseval's identity gives
\[
\|\nabla v\|_{L^2}^2=\|\omega\|_{L^2}^2,
\qquad
\|\nabla^2v\|_{L^2}^2=\|\nabla\omega\|_{L^2}^2,
\]
where the second-derivative norm sums over all ordered
second derivatives. Integrating
\eqref{eq:continuation-vorticity} and applying
Gronwall's inequality to the velocity estimate gives
\begin{equation}\label{eq:continuation-velocity-energy}
	\sup_{0\le t<T_{\max}}
	\bigl(\|v(t)\|_{H^1}^2+\|\omega(t)\|_{L^2}^2\bigr)
	+\int_0^{T_{\max}}\|v(t)\|_{H^2}^2 dt<\infty.
\end{equation}

\noindent{\bf (b) H\"older estimates.}
Fix $0<\tau<T_{\max}$ and use the coordinate
transformation from part 1. The bounds in
\eqref{eq:continuation-velocity-energy} hold in these
variables as well. The two-dimensional Sobolev
embeddings imply
\[
\begin{aligned}
	\|v\cdot\nabla v\|_{L^2(0,T_{\max};L^6)}
	&\le\|v\|_{L^\infty(0,T_{\max};L^{12})}
	\|\nabla v\|_{L^2(0,T_{\max};L^{12})}\\
	&\le C\|v\|_{L^\infty(0,T_{\max};H^1)}
	\|v\|_{L^2(0,T_{\max};H^2)}<\infty.
\end{aligned}
\]
Also,
$\|n(t)\|_{L^6}^6\le
\|n(t)\|_{L^\infty}^4\|n(t)\|_{L^2}^2$.
Hence the right-hand side of the velocity equation
after eliminating the pressure belongs to
$L^2(0,T_{\max};L^6)$.
The evolution operator for its transformed heat part
has Fourier symbol
\[
\exp\left[-\frac1A\int_s^t
\bigl(k^2+(\eta-\sigma k)^2\bigr)d\sigma\right].
\]
By \eqref{eq:local-transformed-ellipticity}, its kernel
is a periodized Gaussian with covariance comparable
to $(t-s)I$. For $0<\gamma<\frac23$, its
$L^6$ to $C^\gamma$ norm is bounded by
$C(t-s)^{-\frac16-\frac{\gamma}{2}}$. Noting that
$\int_0^{T_{\max}}r^{-\frac13-\gamma}dr<\infty$,
Duhamel's formula and Cauchy--Schwarz give a uniform
spatial $C^\gamma$ bound for $v$ on
$[\frac{\tau}{4},T_{\max})$.

For the difference of these evolution operators at
final times $t+h$ and $t$, with the same initial time
$s<t$, Gaussian differentiation gives the
$L^6$ to $L^\infty$ bound
$C\min\{(t-s)^{-\frac16}, h(t-s)^{-\frac76}\}$.
The first bound follows by estimating the operators
separately; the second follows by integrating their
final-time derivatives over an interval of length $h$.
Splitting the Duhamel integral at $t-s=h$ and including
the integral over $[t,t+h]$ gives
\[
\begin{aligned}
	\|v(t+h)-v(t)\|_{L^\infty}
	\le C_\tau h
	+C\left[
	\int_0^h r^{-\frac13}dr
	+h^2\int_h^{T_{\max}}r^{-\frac73}dr
	\right]^{\frac12}\le C_\tau h^{\frac13}
\end{aligned}
\]
for $t\ge\frac{\tau}{4}$, $t+h<T_{\max}$ and
$0<h\le\min\{1,T_{\max}\}$. The initial heat evolution
is bounded by the first term. Larger time increments
are controlled by the uniform supremum bound.
Thus $v$ is uniformly bounded in
$C^{\gamma,\frac{\gamma}{2}}$ on this interval.
The transformed full vorticity satisfies
\[
\partial_t\omega-\frac1A
\bigl[\partial_x^2+(\partial_y-t\partial_x)^2\bigr]\omega
+\frac1A v\cdot\nabla\omega
=-\frac1A\partial_xn.
\]
Its diffusion is uniformly elliptic, its drift is
bounded on $[\frac{\tau}{4},T_{\max})$, and its divergence
forcing has the bounded vector $(-\frac{n}{A},0)$.
The scalar local boundedness and H\"older estimates
in \cite{LSU}, together with the $L^2$ bound in
\eqref{eq:continuation-velocity-energy}, therefore
give boundedness and H\"older continuity of $\omega$
on $[\frac{\tau}{2},T_{\max})$.
The density diffusion matrix is
\[
\frac1A(1+n)^{\alpha-1}[1+(1-\alpha)n]
\begin{pmatrix}1+t^2&-t\\-t&1\end{pmatrix}.
\]
It is bounded and uniformly elliptic because $n$ is
bounded and nonnegative. Its remaining divergence
terms are bounded, since $n$ and $v$ are bounded and
$\|\nabla c(t)\|_{L^\infty}\le C\|n(t)\|_{L^\infty}$.
The same scalar H\"older estimates apply to $n$.
After reducing the exponent, we obtain some
$0<\beta<\gamma$ such that
\begin{equation}\label{eq:continuation-holder-bounds}
	\|n\|_{C^{\beta,\frac{\beta}{2}}}
	+\|\omega\|_{C^{\beta,\frac{\beta}{2}}}
	+\|v\|_{C^{\beta,\frac{\beta}{2}}}<\infty
\end{equation}
on $[\frac{\tau}{2},T_{\max})\times\mathbb T\times\mathbb R$.
These estimates are uniform in space and in
$T<T_{\max}$: they are applied on backward cylinders
of a fixed radius whose bottom times exceed $\frac{\tau}{4}$.

\noindent{\bf (c) Bounds for the first derivatives.}
In the transformed variables, the chemical equation
has the representation
\[
c(t)=\int_0^\infty {\rm e}^{-r}
{\rm e}^{r[\partial_x^2+(\partial_y-t\partial_x)^2]}
n(t)dr.
\]
The gradient heat kernel has $L^1$ norm at most
$Cr^{-\frac12}$, so spatial differences give
$\|\nabla c(t)\|_{C^\beta}\le C\|n(t)\|_{C^\beta}$.
For time differences, use
\[
\partial_t
{\rm e}^{r[\partial_x^2+(\partial_y-t\partial_x)^2]}
=-2r\partial_x(\partial_y-t\partial_x)
{\rm e}^{r[\partial_x^2+(\partial_y-t\partial_x)^2]}.
\]
After one further spatial derivative, the kernel
has $L^1$ norm at most $Crr^{-\frac32}=Cr^{-\frac12}$.
Since $\int_0^\infty {\rm e}^{-r}r^{-\frac12}dr<\infty$,
both the gradient screened kernel and its time
derivative have uniformly bounded $L^1$ norms.
Consequently,
\[
\|\nabla c(t)-\nabla c(s)\|_{L^\infty}
\le C\|n(t)-n(s)\|_{L^\infty}
+C|t-s|\|n(s)\|_{L^\infty}.
\]
Together with \eqref{eq:continuation-holder-bounds},
this proves a uniform $C^{\beta,\frac{\beta}{2}}$ bound for
$\nabla c$ on $[\frac{\tau}{2},T_{\max})$.
The density diffusion matrix and the vector
\[
n(1+n)^\alpha
\begin{pmatrix}1+t^2&-t\\-t&1\end{pmatrix}\nabla c+nv
\]
are now uniformly bounded in $C^{\beta,\frac{\beta}{2}}$.
The interior divergence-form parabolic Schauder
estimate in \cite{LSU}, applied on smaller backward
cylinders, yields
\[
\sup_{\tau\le t<T_{\max}}
\|\nabla n(t)\|_{L^\infty}<\infty.
\]
The coordinate change preserves this spatial
gradient bound.
For the velocity gradient, return to the original
coordinates. At each fixed time the spatial H\"older
bounds for $\omega$ remain uniform. The identities
$v_{\neq}=\nabla^\perp\Delta^{-1}\omega_{\neq}$,
$\partial_yv_{1,0}=\omega_0$,
$v_{2,0}=0$
and the heat representation on non-zero modes give
\[
\begin{aligned}
	\|\nabla v(t)\|_{L^\infty}
	\le\|\omega_0(t)\|_{L^\infty(\mathbb R)}
	+C[\omega_{\neq}(t)]_{C^\beta}
	\int_0^1 r^{-1+\frac{\beta}{2}}dr+C\|\omega_{\neq}(t)\|_{L^2}
	\int_1^\infty {\rm e}^{-\frac{r}{2}}dr.
\end{aligned}
\]
For $0<r\le1$, cancellation of the second derivative
of the heat kernel permits subtraction of the
vorticity value at the evaluation point, giving the
factor $r^{\frac{\beta}{2}}$. For $r>1$, the non-zero $x$
frequencies give exponential decay.
Both integrals are finite. We have therefore proved,
in the original variables,
\begin{equation}\label{eq:continuation-gradient-bounds}
	\sup_{\tau\le t<T_{\max}}
	\bigl(\|\nabla n(t)\|_{L^\infty}
	+\|\nabla v(t)\|_{L^\infty}\bigr)<\infty.
\end{equation}

\medskip
\noindent\textbf{4. Continuation.}
The density bounds and
\eqref{eq:continuation-gradient-bounds} allow us to
apply \eqref{eq:continuation-H3-energy} on $[\tau,T]$,
with $C$ independent of $T<T_{\max}$.
Gronwall's inequality gives
\[
\|n(t)\|_{H^3}^2+\|v(t)\|_{H^3}^2
\le\bigl(\|n(\tau)\|_{H^3}^2+\|v(\tau)\|_{H^3}^2\bigr)
{\rm e}^{C(t-\tau)}.
\]
Integrating \eqref{eq:continuation-H3-energy} then yields
\[
\sup_{\tau\le t<T_{\max}}
\bigl(\|n(t)\|_{H^3}^2+\|v(t)\|_{H^3}^2\bigr)
+\int_\tau^{T_{\max}}
\bigl(\|n(t)\|_{H^4}^2+\|v(t)\|_{H^4}^2\bigr)dt
<\infty.
\]
In the transformed variables, the equations imply
$\partial_tn,\partial_tv\in L^2(\tau,T_{\max};H^2)$.
Applying the time-continuity theorem used in
\eqref{eq:local-H3-time-continuity} on
$[\tau,T_{\max}]$ gives strong $H^3$ limits at
$T_{\max}$. The inverse coordinate and velocity
transformations are uniformly bounded and strongly
continuous on $H^3$ for finite times, so the original
fields have the same terminal regularity.
The terminal velocity is divergence free, and the
terminal density is nonnegative. Fatou's lemma gives
\[
\int_{\mathbb T\times\mathbb R}
n(T_{\max},x,y)dxdy\le M.
\]
Thus the terminal data satisfy the hypotheses of
the local existence result in parts 1--2.
Restarting the original system at $T_{\max}$ extends
the solution beyond $T_{\max}$, contradicting
maximality. This proves the continuation assertion.
\end{proof}
\section*{Declarations}
\subsection*{Acknowledgments}
The author thanks Professor Wendong Wang for his guidance
and helpful discussions.
\subsection*{Competing interests}
The author has no relevant financial or non-financial interests to disclose.
\subsection*{Data availability}
Data sharing is not applicable to this article as no datasets were generated or analysed during the current study.

\end{document}